\documentclass{amsart}

\usepackage{amsmath}
\usepackage{amssymb}

\def\derpar#1#2{\frac{\partial#1}{\partial#2}}
\def\pa{\partial}
\def\var{\varepsilon}
\def\R{{\Bbb R}}

\newtheorem{Thm}{Theorem}
\newtheorem{Prop}[Thm]{Proposition}
\newtheorem{Cor}{Corollary}[Thm]
\newtheorem{Lem}{Lemma}
\def\N{{\Bbb N}}

\begin{document}

\title[Entropy dissipation estimates for the Landau equation]{Entropy dissipation estimates\\ for the Landau equation\\ in the Coulomb case and applications}

\author{L. Desvillettes}

\address{CMLA, ENS Cachan,  CNRS\\
  61, avenue du Pr\'esident Wilson, F-94230 Cachan, France\\
  E-mail: desville@cmla.ens-cachan.fr}

\begin{abstract}
We present in this paper an estimate which bounds from below the entropy dissipation
$D(f)$ of the  Landau operator with Coulomb interaction by a weighted $H^1$ norm 
of the square root of $f$. As a consequence, we get a weighted $L^1_t(L^3_v)$ estimate for the 
 solutions of the spatially homogeneous Landau equation with Coulomb interaction, and the
propagation of $L^1$ moments of any order for this equation. We also present an application
of our estimate to the Landau equation with (moderately) soft potentials, providing thus a new proof of some recent
results of \cite{kcw}.

\end{abstract}

\keywords{Landau equation, Landau operator, entropy dissipation, degenerate diffusion, Coulomb interaction}

\maketitle




\section{Introduction and main result} \label{intro}

\subsection{Description of the Landau equation}


We recall the spatially homogeneous Landau equation of plasma theory (cf. \cite{chapman,lifschitz}),
\begin{equation} \label{Landau}
\derpar{f}{t}  (t,v) = Q(f,f) (t,v), \qquad v \in {\Bbb R}^N,
\quad t \geq 0,  
\end{equation}
where $f:= f(t,v)$ is a nonnegative function representing the number of particles which at time $t$ move with velocity $v$, and  $Q$ is a nonlinear quadratic operator modeling the collisions between those particles.
It acts on the variable $v$ only, and writes
\begin{equation} \label{operator}
Q(f,f)(v) = \sum_{i=1}^N \derpar{}{v_i}   \left \{ \int_{{\Bbb R}^N} \, \sum_{j=1}^N a_{ij}
(v-w)
\left ( f(w) \,\derpar{f}{v_j}(v) - f(v) \,
\derpar{f}{v_{j}}(w) \right )\,  dw\right \} .
\end{equation}
We also introduce the (nonnegative) initial datum $f_{in} : \R^N \to \R_+$:
\begin{equation} \label{donnee_in}
\forall v \in {\Bbb R}^N, \qquad f(0,v) = f_{in}(v).  
\end{equation}
\par
Here, $(a_{ij}(z))_{ij}$ ($z \in {\Bbb R}^N$) 
is given by 
\begin{equation} \label{annee}
a_{ij}(z) = \Pi_{ij}(z) \,\psi(|z|), 
\end{equation}
where $\psi$ is a nonnegative function,  and
\begin{equation}
\Pi_{ij}(z) = \delta_{ij} - \frac{z_i z_j}{|z|^2} 
\end{equation}
is the $i,j$-component of the orthogonal projection $\Pi$
 onto $z^\bot := \{ y\, / \; y\cdot z =0 \}.$
\par
It is customary to define the following functions:
\begin{equation} \label{coeff1}
b_i(z) = \sum_{j=1}^N \partial_j a_{ij} (z) = -(N-1)\, \frac{z_i}{|z|^2}\, \psi(|z|),
\end{equation}
 \begin{equation} \label{coeff2}
c(z) =  \sum_{i=1}^N \sum_{j=1}^N  \partial_{ij} a_{ij} (z)
 =  -(N-1) \,\left( (N-2)\, \frac{\psi(|z|)}{|z|^2} +  \frac{\psi'(|z|)}{|z|}\right) . 
\end{equation}
The computation above must be adapted when $\psi$ has strong singularities (for example, in the important case when $N=3$ and $\psi(|z|) = |z|^{-1}$, corresponding to the Coulomb interaction, one finds $c(z) = -8\pi\, \delta_0$).
\medskip

Using those functions, the Landau operator can be written (at the formal level)
as a (conservative or non conservative)
nonlinear (quadratic) parabolic equation with nonlocal coefficients, indeed
\begin{equation} \label{parab1}
  Q(f,f)= \sum_{i=1}^N \frac{\pa}{\pa v_i} \bigg(  \sum_{j=1}^N  (a_{ij}  * f) \, \frac{\pa f}{\pa v_j}
- (b_i * f)\, f \bigg),
\end{equation}
and
\begin{equation} \label{parab2}
  Q(f,f)= \sum_{i=1}^N  \sum_{j=1}^N  (a_{ij}  * f)\, \frac{\pa^2 f}{\pa {v_i}\pa{v_j}} 
- (c * f)\, f . 
\end{equation}

We also observe that at the formal level (that is, when both $f$ and $\varphi$ are smooth functions having a reasonable 
behavior at infinity), the following weak version of the
 Landau operator can be defined:

\begin{equation} \label{weakformop}
 \int_{\R^N} Q(f,f)(v) \,\varphi(v)\, dv 
\end{equation}
$$ = \frac{1}{2} \sum_{i=1}^N  \sum_{j=1}^N \iint_{\R^N\times\R^N}
 f(v)\, f(w)\, a_{ij}(v-w)\, 
\Bigl ( \partial_{ij} \varphi(v) + (\partial_{ij}\varphi)(w) \Bigr ) \, dv dw $$
$$ \qquad\qquad {} +   \sum_{i=1}^N  \iint_{\R^N\times\R^N} f(v)\, f(w)\, b_i(v-w)\, 
\Bigl ( \partial_i \varphi(v) - (\partial_i \varphi)(w) \Bigr )\, dv dw. $$

Using the test functions $\varphi(v) = 1$, $v_i$ (for $i=1,..,N$), $\frac{|v|^2}2$,
 we see that (still at the formal level), the Landau operator
conserves mass, ($i$-th component of the) momentum and kinetic energy, that is:
\begin{equation} \label{conserv}
 \int_{\R^N} Q(f,f)(v) \, \left(\begin{array}{c} 1\\ v_i\\ |v|^2/2 \end{array} \right)\, dv = 0. 
\end{equation}
We also get (at the formal level) 
the formula for the entropy dissipation $D_{\psi} := D_{\psi}(f)$ (defined on functions $f$ from $\R^N$ to $\R_+$) by considering $\varphi(v) = \ln f(v)$:
\begin{equation}\label{mded} 
D_{\psi}(f)  := -\int_{{\Bbb R}^N}  Q(f,f)(v) \ln f(v) \, dv 
\end{equation}
\begin{multline*}=\frac12\, \sum_{i=1}^N \sum_{j=1}^N \iint_{{\Bbb R}^N\times {\Bbb R}^N} 
f(v)\,f(w)\,
\psi(|v-w|) \, \Pi_{ij}(v-w)\, \left ( \frac{\partial_i f}{f}(v) -
\frac{\partial_i f}{f}(w) \right ) \\
\qquad \qquad  \left( \frac{\partial_j f}{f}(v) -
\frac{\partial_j f}{f} (w) \right )\, dv dw 
\end{multline*}
$$ = 2 \iint_{{\Bbb R}^N\times {\Bbb R}^N}  \psi(|v-w|) \, 
\left| \Pi\, [ \,(\nabla_v - \nabla_w) \sqrt{f(v)\,f(w)} \,]\, \right|^2\, dvdw \ge 0. $$
Note that this last formula (as already noticed in \cite{vill:new:97}, 
but with $\psi$ inside the derivatives) implies that $D_{\psi}(f)$
is defined (but may be infinite, in particular if the distribution $\Pi\,(\nabla_v - \nabla_w) \sqrt{f(v)\,f(w)}$ 
does not lie in $L^2_{loc}(\R^N\times\R^N)$) as soon as $f \in L^1_{loc}(\R^N)$.
\medskip

As a consequence of formulas (\ref{conserv}) and (\ref{mded}), 
the solutions of the Landau equation (\ref{Landau}), (\ref{donnee_in})
satisfy (at least formally)
the conservation of mass, momentum and energy, that is
\begin{equation}\label{consmm} 
 \int_{{\Bbb R}^N}
f(t,v) \, \left(\begin{array}{c} 1\\ v_i\\ |v|^2/2 \end{array} \right)\, dv = \int_{{\Bbb R}^N} f_{in}(v)\,
\left(\begin{array}{c} 1\\ v_i\\ |v|^2/2 \end{array} \right)\, dv,
\end{equation}
and the entropy dissipation identity 
\begin{equation}\label{ED}
\frac{d}{dt} H(f(t,\cdot)) 
 = - D_{\psi}(f(t,\cdot))\le 0, 
\end{equation}
where $H := H(f) $ is the entropy functional (defined on functions from $\R^N$ to $\R_+$)
\begin{equation} \label{ent_diff}
H(f) := \int_{{\Bbb R}^N} f(v) \ln f(v)\, dv,
\end{equation}
and $D_{\psi}$ is the entropy dissipation functional defined in (\ref{mded}).
\par
Note that identities (\ref{consmm}) and (\ref{ED}) 
 naturally furnish an {\it{a priori}} estimate: indeed integrating (\ref{ED})
on the time interval $[0,T]$, one can write
$$  H(f(T,\cdot)) 
  + \int_0^T D_{\psi}(f(t,\cdot))\, dt = H(f(0,\cdot)). $$
  Then, 
  $$ \int_{\R^N} f(t,v)\, |\ln f(t,v)|\, dv = \int_{\R^N} f(t,v)\, \ln f(t,v)\, dv$$
$$  +\, 2\,  \int_{\R^N} f(t,v)\, (-\ln f(t,v))\, 1_{\{ f(t,v) \le 1\} } \,  dv $$
  $$ = H(f(t,\cdot)) + 2\,  \int_{\R^N} f(t,v)\, (-\ln f(t,v))\, 1_{\{f(t,v) \le e^{-1 - |v|^2} \} } \,  dv $$
  $$ + \,2\,  \int_{\R^N} f(t,v)\, (-\ln f(t,v))\, 1_{\{ 1\ge f(t,v) \ge e^{-1- |v|^2} \} } \,  dv $$
  $$ \le H(f(t,\cdot)) +  2\,  \int_{\R^N}  e^{-1 - |v|^2}\, (1+|v|^2)\, dv +  2\,\int_{\R^N} f(0,v)\, (1+|v|^2) \,dv. $$
Then, for any initial datum such that the mass, energy and entropy are finite, that is
\begin{equation}\label{mee}
{\mathcal{M}}_{in} := \int_{v\in\R^N} f_{in}(v)\, \bigg( 1 + \frac{|v|^2}2 + |\ln (f_{in}(v))| \bigg) \, dv < + \infty,
\end{equation}
we see that a solution of (\ref{Landau}), (\ref{donnee_in}) satisfies at the formal level the following
 {\it{a priori}} estimate:
 \begin{equation}\label{meeap}
 \sup_{t\in [0,T]} \int_{v\in\R^N} f(t,v)\, \bigg( 1 + \frac{|v|^2}2 + |\ln f(t,v)| \bigg) \, dv
 \end{equation}
$$ + \int_0^T D_{\psi}(f(t,\cdot))\, dt \le C(T,N, {\mathcal{M}}_{in}), $$
where the constant $C(T,N, {\mathcal{M}}_{in})$ only depends on $T,N$ and $ {\mathcal{M}}_{in}$.
\bigskip

The most important (that is, physically relevant) function $\psi$ appearing in operator (\ref{Landau})
is $\psi(z) = |z|^{-1}$ in dimension $3$. It corresponds to the case of charged particles (moving according to Coulomb interaction), cf. \cite{lifschitz}, and it also naturally appears in the so-called weak coupling asymptotics of Boltzmann equation (cf. \cite{boby} and the older reference \cite{bogol}). Our paper is mainly devoted to the study of
this case (systematically called in the sequel the ``Coulomb case''). In the Coulomb case, $b_i(z) = - 2\, z_i\, |z|^{-3}$, and $c(z) = - 8\pi\,\delta_0(z)$.
\bigskip

It is however also interesting, at least from the mathematical viewpoint, to consider more general 
functions $\psi$ and dimensions $N \in \N$ ($N\ge 2$). We refer for example to \cite{desv_TTSP} to see 
how the Landau kernel with arbitrary $\psi$ can be obtained from the Boltzmann kernel (with arbitrary cross section) through a  scaling involving the concept of grazing collisions.
\par 
We use in this paper a terminology very close to that of \cite{kcw}. When $N=3$ (and sometimes for all $N \in \N$, $N\ge 2$), if $\psi$ is given by a power law, we say that 
\begin{equation}\label{pote}
\psi(|z|)= |z|^{\gamma+2}
\end{equation}
  is coming out of hard potentials when $\gamma \in ]0,1]$, Maxwellian molecules when $\gamma =0$, (moderately) soft potentials when $\gamma \in [-2,0]$, and very soft potentials 
when $\gamma \in ]-4,-2[$. Note that the Coulomb case falls within the category of very soft potentials such as it is defined here, and that it does not constitute a limiting case of this category. 
\par
As we shall see, one of the main results of this paper holds as long as  $\gamma > -4$; the limiting case corresponds 
thus to $\gamma = -4$. This provides some justification for the use of the terminology ``very soft potentials'' in the
case when $\gamma \in ]-4,-2[$, rather than for another interval whose right-hand side is $-2$. Note that $\gamma=-4$ already appeared as a limiting parameter in the context of the Boltzmann equation, cf. applications in section 7 of \cite{advw}.

 \subsection{Existing results on the spatially homogeneous Landau equation and new results in the Coulomb case}
\medskip

Our primary goal in this paper is to improve the existing theory on the spatially homogeneous Landau equation
in the Coulomb case. We however first recall here what is known at the mathematical level in ``easier'' cases, and only
describe at the end what is known in the specific case of the Coulomb interaction (in 3D).
\smallskip

 The case of hard potentials, (that is, $\psi$ given by (\ref{pote}) with $\gamma \in ]0,1]$), has been thoroughly treated. We refer to \cite{desvvill:I:97} for the existence of 
$C^{\infty}$ solutions, and to \cite{chenlixu1}, \cite{chenlixu2} and \cite{mori} for a study of smoothness beyond $C^{\infty}$.
The solutions are moreover unique, and rapidly decaying (w.r.t. the $v$ variable, and when $t>0$) 
whenever the initial datum lies in a suitably weighted $L^2$ space.
Finally, Maxwellian lower bounds are also known to hold in this case (cf. \cite{desvvill:I:97}).
\smallskip

The special case $\gamma = 0$ in (\ref{pote}) (that is, Maxwellian molecules) is also quite well-known (cf. \cite{vill:lm:96}): smoothness also immediately appears when 
the initial datum has 
finite mass and energy, and uniqueness of a classical solution holds for such an initial datum. Polynomial moments (with respect to the $v$ variable) are however only propagated (not created like when $\gamma>0$). Maxwellian lower bounds exist as in the case of hard potentials.
A remarkable feature of Maxwellian molecules is that the Fisher information 
is
decaying, so that $H$ is not the only known Lyapunov functional in that case.
\smallskip

The theory is less developed in the case of (moderately) soft potentials (that is, $\psi$ given by (\ref{pote}) with $\gamma \in [-2,0[$).
 The best result as far as existence/smoothness is concerned can be found in the very recent paper
\cite{kcw} (cf. also \cite{alex_new} for related estimates, and \cite{guerin} for measure solutions). There,
 propagation of $L^p$ norms is proven, which implies existence and uniqueness
of suitably defined solutions thanks to theorems proven in \cite{fournier_guerin}. Note that the limiting case $\gamma =-2$, which requires a lot of care, is  
 treated in \cite{kcw}. We recover in Appendix 2 of our paper some of the results of \cite{kcw}, with a different proof.
\smallskip

In all the previous cases, the boundedness of the mass, energy and entropy of $f$ 
 is sufficient to give a sense to 
weak solutions (that is, solutions defined by using the weak form (\ref{weakformop}) of the operator). Then, 
extra smoothness estimates are based on the parabolic forms (\ref{parab1}), (\ref{parab2}) of the operator. 
More precisely, coercivity estimates are extracted from the term involving $a_{ij} * f$ in those formula. They take the form (cf. \cite{alex_new}, Prop. 2.1. p.4):
\begin{equation}\label{coerc}
 \sum_{i=1}^N\sum_{j=1}^N (a_{ij} * f)(v) \, \xi_i\, \xi_j \ge C\,(1+|v|^2)^{\gamma/2} \, |\xi|^2, 
 \end{equation}
where $C$ only depends on the mass, energy and entropy of $f$.
Estimate (\ref{coerc}) can be efficient only if the terms 
involving $b_i * f$ in the (conservative) parabolic equation (\ref{parab1}) can be controlled. Such a control is rather easily obtained  in the case of hard potentials 
or Maxwellian molecules ($\gamma \in [0,1]$ in (\ref{pote})), but becomes much more difficult to get when (moderately) soft potentials ($\gamma \in [-2,0[$ in (\ref{pote})) are concerned (cf. \cite{kcw} and Appendix 2 of the present work for more details).
\smallskip

 For very soft potentials ($\gamma \in ]-4,-2[$ when $\psi$ is given by (\ref{pote})), several new difficulties appear. 
\par
First, the mere boundedness
 of mass, energy and entropy of $f$
 does not seem sufficient any more for defining  weak solutions thanks to (\ref{weakformop}).
 Indeed, remembering that
$a_{ij}(v-w)$ behaves like $|v-w|^{\gamma +2}$ (when $|v-w|$ is close to $0$), terms appearing in 
the weak form (\ref{weakformop}) of the operator (written here in 3D, but the same difficulty appears in all dimensions)
such as
 \begin{equation}\label{termad}
\iint_{\R^3\times\R^3}
 f(v)\, f(w)\, a_{ij}(v-w)\, 
\Bigl ( \partial_{ij} \varphi(v) + (\partial_{ij}\varphi)(w) \Bigr ) \, dv dw 
\end{equation}
are not automatically  defined when $f$ belongs to a (weighted) $L^1(\R^3)$ (or $L\,\ln L(\R^3)$) space. Note that they are however defined if $f$ belongs to some $L^p(\R^3)$ for $p>1$
large enough. Looking for the best (smallest) possible $p$ which enables the definition of (\ref{termad}), for $f$ which do not depend on the time variable $t$,
 gives an insight
of the specific difficulties of very soft potentials. Observing that when $\gamma < -2$,  $|\cdot|^{\gamma +2}\,1_{|\cdot|\le 1}$ lies in 
$L^{\frac{3}{-\gamma -2} -\var}$ for all (small)
 $\var>0$, we see that if $f\in L^p(\R^3)$, then $f *  [a_{ij}\,1_{|\cdot|\le 1}] \in L^r(\R^3)$,
with $\frac1r > \frac1p + \frac{-\gamma -2}{3} - 1$,
 thanks to Young's inequality for convolutions.
 It is then possible to define (\ref{termad})
if $\frac1r + \frac1p < 1$, that is $\frac2p + \frac{-\gamma -2}{3} <2$, or $p>\frac6{8+\gamma}$.
This computation is in fact too optimistic because of the dependence of $f$ w.r.t. time, but it 
already shows that some $L^p(\R^3)$ (for $p>1$) estimate is required in order 
to define weak solutions.
\par
Secondly, the attempts to use the coercivity estimates such as (\ref{coerc}) 
have not been able yet to provide $L^p$ estimates with $p>1$ in the case of very soft potentials ($\gamma \in ]-4,-2[$ in (\ref{pote})), because the term involving $b_i *f$ in (\ref{parab1}) cannot be controlled. This is all
the more frustrating since such $L^p$ estimates exist when $\gamma \ge - 2$ (cf. \cite{kcw}), but  
are not needed to define weak solutions in that case (as already mentioned, an $L^p$ estimate  with $p>1$ is
needed in order to define a weak solution only in the case of very soft potentials).
\smallskip

Note nevertheless that (in the spatially homogeneous as well as in the spatially inhomogeneous context),
it is possible (in the Coulomb case) to build a theory of local (in time) solutions, or of global solutions with small initial data, cf. \cite{arsenev}, \cite{fournier}, \cite{guo} and \cite{alex_new}.
\smallskip

As far as global in time solutions with general (not small) initial data are concerned,
a significant step forward was achieved in \cite{vill:new:97}, where it was observed that the $L^p$ estimate
(with, say, $p$ larger than $\frac6{8+\gamma}$) 
was not really necessary to define solutions, if one allows to use the 
boundedness (in $L^1$ w.r.t. time) of the entropy dissipation $D_{\psi}(f)$ (and not only the mass, energy and entropy appearing in the {\it{a priori}} estimate (\ref{meeap})) to 
bound terms like (\ref{termad}). Solutions
defined in this way are called ``H-solutions'' in order to distinguish them from weak solutions which need
an $L^p$ estimate
(with $p>1$ large enough) for $f$.
\smallskip

Another progress was made in \cite{av} (see formula (70) p. 30, Lemma 13 p. 35, and Remark p. 36), when renormalized solutions of the spatially 
inhomogeneous Landau equation in the Coulomb case were shown to satisfy bounds of the type 
  \begin{equation}\label{bizbiz}
 \int_0^T \iint_{\R^{2N}} \sum_{i,j} (a_{ij} *_v f)(v)\, \pa_{v_i} \left(\frac{f}{1+f}\right)
 \, \pa_{v_j} \left(\frac{f}{1+f}\right) \, dv dx dt < +\infty ,
 \end{equation}
 which imply that the terms appearing in the definition of renormalized solutions 
 (of the spatially 
inhomogeneous Landau equation in the Coulomb case) are well defined.
 \smallskip

 This result is related to the estimates on the non cutoff Boltzmann operator proven
 in \cite{advw} (cf. also the earlier works \cite{lions}, 
\cite{vill:smooth} and \cite{alex_entropy}), which belong to the class of entropy dissipation estimates. This means that the entropy dissipation functional related to a PDE is bounded below by a functional controlling somehow the smoothness
of the function (appearing in the entropy dissipation functional). 
Note that the results of \cite{advw} were very significantly improved in
\cite{gs1}, \cite{gs2}, \cite{gs3} thanks to the introduction of an anisotropic distance specifically adapted to the geometry of collisions in the Boltzmann equation.
\smallskip

It is in fact proven in
\cite{advw} that the entropy dissipation of the Boltzmann operator without cutoff controls (with
the help of the mass and energy)
$\sqrt{f}$ in $H^s_{loc}$, where $s$ is related to the strength of the angular singularity 
in the cross section of the Boltzmann equation. By suitably scaling 
this cross section (cf. \cite{advw} or \cite{av}), it is
possible to recover the Landau equation, thus suggesting that in the limit, the 
Landau entropy dissipation (with, say, $\psi$ given by (\ref{pote}) and $\gamma \in ]-4,1]$), controls  $\sqrt{f}$ in $H^1_{loc}$. Such a result is indeed quoted  
in \cite{villani_hb} (section 5.1.3.2, p. 170), and may be obtained either starting from 
the results of \cite{advw}, or from those of \cite{DV2}.
Using then a Sobolev embedding, the {\it{a priori}} estimate (\ref{meeap}) implies
that (a solution of the spatially homogeneous Landau equation) $f$ lies in $L^1([0,T]; L^3_{loc}(\R^N))$, so that the weak formulation (\ref{weakformop}) makes sense. 
\smallskip

In this paper, we show that in fact, in the Coulomb case, it is possible to go beyond a local setting,
and get a
global weighted (presumably with the optimal weight) estimate in
$L^1([0,T]; L^3(\R^3))$ for $f$,  associated to a global weighted (presumably with the optimal weight) bound of $\sqrt{f}$  in $L^2([0,T]; H^1(\R^3))$.

 This estimate is a consequence of our main theorem, which is an entropy dissipation estimate
enabling the control of the weighted $H^1(\R^3)$ norm of $\sqrt{f}$ by the Landau  entropy dissipation $D_{|\cdot|^{-1}}$ of the Landau operator in the Coulomb case.
\medskip

We introduce the following definition for weighted $L^p$ spaces:
for all $l \in\R$, $p\in [1, +\infty]$, the weighted $L^p$ spaces and norms are thus defined: 
$$ ||f||_{L^p_l(\R^N)} := || (v \mapsto (1+|v|^2)^{l/2} \, f||_{L^p(\R^N)}, $$
and $ L^p_l(\R^N) = \{f: \R^N \to \R, \, ||f||_{L^p_l(\R^N)} < +\infty \}$.
\medskip

Then, our main theorem writes
\smallskip


\begin{Thm}\label{princ2}
Let  $f := f(v)\ge 0$, belonging to $L^1_2(\R^3)$,
be such that 
 $\int f\, |\ln f|\, dv \le \bar{H}$, for some $\bar{H}>0$.
\smallskip

Then, there exists a constant $C:= C(\int f\, dv, \int f\,v\,dv, \int f\,|v|^2/2\,dv,\bar{H})>0$ which (explicitly) depends only on the mass $\int f\, dv$, the momentum $\int f\,v\,dv$, 
the energy $\int f\,|v|^2/2\,dv$ and (an upper bound on the) entropy $\bar{H}$,  such that
$$ \int_{\R^3} |\nabla \sqrt{f(v)}|^2\, (1+ |v|^2)^{- 3/2} \, dv  \le  C\,(1+ D_{x \mapsto |x|^{-1}}(f)).$$
 If moreover $f$ is radially symmetric and satisfies the normalization
 $\int f(v)\, dv =1$, $\int f(v)\, v\, dv =0$, $\int f |v|^2\, dv =3$, then the constant can be explicitly  estimated with numbers which are not too huge:
$$  \int_{\R^3} |\nabla \sqrt{f(v)}|^2\, (1+ |v|^2)^{- 3/2} \, dv  \le   108$$
$$ \times\, 13^{3/2} \, \left(\frac{16\pi}3 \right)^{4/3} \,
\exp\left(\frac{16}3\, \bar{H}\right)\, \left(2+ \frac{128}3\, D_{x \mapsto |x|^{-1}}(f) \right).$$
\end{Thm}

This theorem is proven in section \ref{sec2}, as a corollary of a more general theorem presented there (Theorem \ref{princ1}), in which the dimension 
is arbitrary, and the hypothesis on $\psi$ is relaxed: we assume indeed only that $\psi$ is bounded below by 
a (negative) power law at infinity. The discussion on the hypothesis of this theorem, as well as a presentation
of the main arguments of the proof (and the link of this proof with other works) is postponed to section \ref{sec2}. 
\medskip

As a corollary of Theorem \ref{princ2} (and as previously announced), we recover that H-solutions of the Landau equation in the Coulomb case
(defined on $[0,T]$) belong to $L^1([0,T], L^3_{-3}(\R^3))$, and satisfy the estimate
$$  \int_0^T \int_{\R^3} |\nabla \sqrt{f(t,v)}|^2\, (1+ |v|^2)^{- 3/2} \, dv dt  < +\infty. $$
 They are therefore ``usual'' weak 
solutions (of the Landau equation in the Coulomb case). 
This last result (the fact that $H$-solutions are usual weak solutions) was in fact already suggested in \cite{advw}, where a close computation was performed 
for the Boltzmann equation without cutoff and very soft potentials (cf. second application of section 7 in \cite{advw}). The bound suggested there (and in \cite{villani_hb}, section 5.1.3.2, p. 170) is however only local (that is,
$ f\in L^1([0,T], L^3_{loc}(\R^3))$).
\medskip

More precisely, we prove the
\medskip


\begin{Cor} \label{weakso}
Let $T>0$,  $f_{in} := f_{in}(v) \in L^1_2({\Bbb R}^3) \cap L\,\ln L(\R^3)$ be nonnegative,
and $f:= f(t,v) \in L^{\infty}([0,T];  L^1_2(\R^3) )$ be a (nonnegative) H-solution of the Landau equation in the Coulomb case ($\gamma = -3$ and $N=3$) with initial datum $f_{in}$ (such a solution exists thanks to \cite{vill:new:97}).
\par
Then, $f\in L^1([0,T], L^3_{-3}(\R^3))$, 
%
and for all $\displaystyle \varphi := \varphi(t,v) \in
C^2_c([0,T[ \times {\Bbb R}^3_v)$,
\begin{equation} \label{star}
  - \int_{\R^3} f_{in}(v)\, \varphi(0,v)\, dv - \int_0^{T} \int_{\R^3} 
f(t,v)\, \partial_t \varphi(t,v)\, dv dt
\end{equation}
$$ = \frac{1}{2} \sum_{i=1}^3  \sum_{j=1}^3 \int_0^T \iint_{\R^3\times\R^3}
 f(t,v)\, f(t,w)\, a_{ij}(v-w)\, 
\Bigl ( \partial_{ij} \varphi(t,v) + (\partial_{ij}\varphi)(t,w) \Bigr ) \, dv dw\, dt $$
$$  \qquad\qquad {} +   \sum_{i=1}^3  \int_0^T\iint_{\R^3\times\R^3} f(t,v)\, f(t,w)\, b_i(v-w)\, 
\Bigl ( \partial_i \varphi(t,v) - (\partial_i \varphi)(t,w) \Bigr )\, dv dw\, dt. $$
\bigskip

We underline that all the terms above (remember that $a_{ij}(z) = \Pi_{ij}(z)\, |z|^{-1}$ and
$b_i(z) = -2 \,z_i\,|z|^{-3}$) are defined under the assumption that
$ f\in L^\infty([0,T]; L^1_2({\Bbb R}^3)) \cap
L^1 ([0,T]; L^{3}_{-3}({\Bbb R}^3)))$. This last statement is part of Corollary \ref{weakso}.
\end{Cor}

Note that thanks to an interpolation (see Proposition \ref{interp3}) 
between $L^\infty([0,T]; L^1_2({\Bbb R}^3))$ and
$L^1 ([0,T]; L^{3}_{-3}({\Bbb R}^3)))$, we get for $f:=f(t,v)$ (weak solution of the Landau equation in the Coulomb case)
 an estimate in $L^{1+\delta}([0,T];L^{3 -\var}_{-3}(\R^3))$ for some $\var, \delta >0$,  so that the first estimate proven in 
\cite{kcw} in the case of (moderately) soft potentials still holds for the Coulomb case, up to the
weight (this restriction is removed in next paragraph). More details about
the link between this work and the estimates of \cite{kcw} are presented in Appendix 2, let us only add here that our estimate is uniform in time (whereas in \cite{kcw} it grows polynomially or more than polynomially (when $\gamma =-2$)) in time, but it involves a negative weight in $v$.
Removing this weight thanks to an interpolation with moments in $L^1$ (see paragraph below) 
 leads, as in \cite{kcw}, to a polynomially growing w.r.t. time estimate. 
\medskip

We also indicate that our result does not allow to use the uniqueness theorems of \cite{fournier_guerin} and \cite{fournier}, which state that for $\gamma \in [-3,-2[$, uniqueness holds
(for solutions of the spatially homogeneous Landau equation) when the solution $f$ lies
in $L^{\infty}([0,T]; L^1_2(\R^3)) \cap L^1 ([0,T]; L^{3/(3 + \gamma)}(\R^3))$ 
(with $3/(3 + \gamma)$ replaced by $\infty$ when $\gamma =-3$, that is, in the Coulomb case).
Indeed, $3/(3 + \gamma) > 3$ whenever $\gamma \in [-3,-2[$.
\medskip

Our main estimate (Theorem \ref{princ2}) can also be used in order to study the behavior of weak solutions 
of the Landau equation in the Coulomb case when $|v| \mapsto +\infty$. As very often when collision operators 
are concerned, this behavior is studied through the evolution 
of the ($L^1$) moments of the solution. As expected for (moderately or very) soft potentials, we are able to show that moments of any order are propagated when $H$-solutions (or equivalently, weak solutions) of the Landau equation in the Coulomb case are considered. More precisely, we show the

\begin{Prop} \label{momco}
Let $T>0$, $k\in \R$, $k>2$, and  $f_{in} := f_{in}(v) \in L^1_k({\Bbb R}^3)$ be nonnegative.
We consider  $f:= f(t,v) \in L^{\infty}([0,T]; L^1_2(\R^3))$  a (nonnegative)  H-solution (or equivalently, weak solution) of the Landau equation in the Coulomb case ($\gamma = -3$ in (\ref{pote}) and $N=3$), with initial datum $f_{in}$.
\par
Then, $f\in L^{\infty}([0,T]; L^1_k(\R^3))$ for all $T>0$ (and the corresponding norm can be explicitly estimated in terms
of $T$, the initial mass, momentum, energy, (upper bound of) entropy, and initial $L^1_k({\Bbb R}^3)$ bounds).
\end{Prop}

This proposition is the consequence of a more general
 proposition (Proposition~\ref{momco}) proven in section \ref{momse}, in which the assumptions on the function $\psi$ are relaxed, and in which all dimensions are considered. Once again, the discussion on the assumptions
 of the proposition, and the link with other results on propagation of moments is postponed to section \ref{momse}.
\smallskip

 We only write here this small comment: this result shows that, as expected,
the behavior for large $|v|$ of the solutions of the Landau equation in the Coulomb case
 is similar to what is observed for (moderately) soft potentials: polynomial moments are propagated but not
created (contrary to the case of hard potentials).

  \subsection{Structure of the paper}

We start with the proof of our main theorem (Theorem \ref{princ2} and its more general version Theorem \ref{princ1}), together with various comments on this theorem and its proof, in section
\ref{sec2}. 
\par
Corollary \ref{weakso} is then proven and commented in section \ref{weaksose}.
\par
 Section \ref{momse} is devoted to the proof of the statement on moments (Proposition \ref{momco}).
\par
 Finally, we add some appendices 
for the sake of completeness: standard interpolations are presented in Appendix 1, while some remarks 
on the links between the results presented here and the regularity issues on (moderately) soft 
potentials discussed in \cite{kcw}, are presented in Appendix 2. 

\section{The entropy dissipation controls the (weighted) Fisher information} \label{sec2}

\subsection{Description of our main result in a general setting}

 We first write down a theorem which is slightly more general than Theorem \ref{princ2}, since it concerns
 all dimensions of space, and general functions $\psi$:

\begin{Thm} \label{princ1}
Let  $f := f(v)\ge 0$, belonging to $L^1_2(\R^N)$,
be
such that 
 $\int f\, |\ln f|\, dv \le \bar{H}$, for some $\bar{H}>0$.
Let $\psi$ satisfy
$$ \forall z\in \R^N, \qquad \psi(z) \ge K_3\, \inf(1, |z|^{\gamma_1 + 2}), $$
for some $K_3>0$ and $\gamma_1 \le 0$.
\smallskip

Then, there exists a constant $C:= C(N, \int f\, dv, \int f\,v\,dv, \int f\,|v|^2/2\,dv,\bar{H}, \gamma_1,K_3)>0$ which (explicitly) depends only on the dimension,
 the mass, momentum, energy, (an upper bound of the) entropy
and the parameters of the lower bound on $\psi$:  $\gamma_1$, $K_3$,
 such that
$$  \int_{\R^N} |\nabla \sqrt{f(v)}|^2\, (1+ |v|^2)^{\inf(\gamma_1/2, -1)} \, dv  \le  C
 \, (1+ D_{\psi}(f)),$$
where $D_{\psi}(f)$ is defined in (\ref{mded}).
\end{Thm}

{\bf{Remarks}}: The assumptions on $\psi$ in this theorem are not optimal. Indeed,
 one can check that as long
as $\psi$ is bounded below by any function decaying less rapidly than a Maxwellian function of $v$, the 
same result can be obtained (with a slight modification of the proof given in this paper), up to a change of weight. Our feeling is that even a Maxwellian lower bound is not compulsory (though the proof should then be changed more deeply), if one allows a weight tending to $0$ very quickly at infinity.
\smallskip

It is clear that $D_\psi$ behaves in a monotone way w.r.t $\psi$, so that one could hope to get a better 
estimate when $\psi$ has a singularity (at point $0$). This cannot be seen on the estimate presented here,
in which we use $1$ as a bound from below for $\psi$ at finite distance. Therefore,
 there might be room for improvement of the estimate when $\psi$ is given by a
 power law (and especially in the Coulomb case).
\smallskip

The weight appearing in $\int_{\R^N} |\nabla \sqrt{f(v)}|^2\, (1+ |v|^2)^{\inf(\gamma_1/2, -1)} \, dv$ 
behaves like $|v|^{\gamma_1}$ for large $|v|$ (as soon as $\gamma_1 < -2$), whereas $\psi$ behaves like $|v|^{\gamma_1 +2}$. This loss
of weight is reminiscent of the loss of weight in coercivity estimates (cf. for example \cite{desvvill:I:97}
and \cite{kcw}), which illustrates the degenerate character of the elliptic operator appearing
in the parabolic formulations of the equation. Our feeling is therefore that this weight is somehow optimal.
\bigskip


The proof relies on computations which link the entropy dissipation (of the Landau equation) to quantities like $ \frac{\nabla f}f$. Those computations, which use integrations against 
test functions of the type $f(w)$ or $w\,f(w)$, were used in the study of the long time
behavior of the Landau equation (cf. \cite{DV2}), and are themselves reminiscent of old computations linking the entropy dissipation to the distance to the space of Maxwellian states 
(cf. \cite{desv:entr:89}), where derivations instead of integrations were used.    
\smallskip

An important difference between the proof of Theorems \ref{princ2}, \ref{princ1} and those previous works is that because we want to treat $\psi$ which decay like a negative power law at infinity, and use only the moments of order less than $2$ of $f$ (because only those are {\it{a priori}} propagated by the Landau equation), one needs to change the test functions used
in those previous works, that is 
$f(w)$ and $w\,f(w)$, in functions which decay fast at infinity, like   
$f(w)\, e^{- \delta\, |w|^2}$ and $w\,f(w)\,e^{- \delta\, |w|^2}$, for some
$\delta>0$. This leads to a new 
difficulty related to the fact that quantities like $\int  w\,f(w)\,e^{- \delta\, |w|^2}\, dw$
are not as easy to manipulate as $\int  w\,f(w)\,dw$ (which appear in the assumption of the 
theorem). 
\smallskip

In particular, in the course of the proof, one needs (for some $\lambda>0$) a lower bound on the quantity
$$ \bigg| {\hbox{ Det }} \left( \int_{\R^N} e^{- \lambda w^2}  \left[ \begin{array}{ccc}
 1  & w_j & w_i \\
 w_i  & w_i\,w_j & w_i^2 \\
 w_j & w_j^2& w_i\,w_j\end{array} \right]  \,f(w)\, dw \right) \,\bigg| , $$
provided that $f$ has a given (non zero) mass and energy, a given momentum and an upper bound on its entropy.
\par
The same result with $\lambda =0$ was obtained in \cite{DV2}, thanks to the 
fact that by a change
of functions and variables, one can impose that the momentum is $0$ (that is, $\int f(w)\,w\,dw =0$), and that the matrix of directional temperatures $\left(\int_{\R^N} w_i\,w_j\, f(w)\, dw \right)_{i,j =1,..,N}$ is diagonal. 
\par
The proof of \cite{DV2} seems difficult to modify in order to take into account the extra term $e^{- \lambda w^2}$
in the integrals. We therefore provide a new proof, based on the direct computation of the
determinant and the fact that one can look for a parameter $\lambda>0$ as small as desired.
\par
However, in the special case of radially symmetric $f$, one has  $$\int  w\,f(w)\,e^{- \delta\, |w|^2}\, dw = 0, $$ and the proof of \cite{DV2} can be used after suitable modifications. Theorem
\ref{princ2} is then much easier to prove. We provide therefore this simplified proof in next subsection, before presenting the complete proof of Theorems \ref{princ2} and \ref{princ1}.
\bigskip

 The rest of this section is thus dedicated to the proof of Theorems \ref{princ2} and \ref{princ1}.
A simplified proof of Theorem \ref{princ2} (that is, the
 Coulomb case) is first presented (subsection \ref{sub1}). It holds only when radially symmetric distribution 
 functions $f$ are considered. We also restrict ourselves in this first proof to 
the case when $f$ has a normalized mass, momentum and energy, in order to simplify it further, and
to keep track of numerical constants.
\par
 The complete proof of Theorems \ref{princ2} and \ref{princ1} is then presented in subsection \ref{sub2}.
 \bigskip
 

 \subsection{Proof in a simplified case} \label{sub1}

In this subsection, we first prove Theorem \ref{princ2} when $f$ is supposed to be radially symmetric
and satisfies the normalization
\begin{equation}\label{norma}
\int_{\R^3} f(v)\, dv =1, \qquad \int_{\R^3} f(v)\, v\, dv =0, \qquad \int_{\R^3} f |v|^2\, dv =3.
\end{equation}
We also assume that $f$ has a sufficient (strictly positive) bound from below for the computations below to make sense. 
\bigskip

{\bf{Proof of Theorem \ref{princ2} when $f$ is supposed to be radially symmetric and satisfies the normalization (\ref{norma})}}:
 
We first observe that (for all $x,y\in\R^3$)
$$ y^T\,(|x|^2 \, Id- x \otimes x)\,y = \frac12 \, \sum_{i,j=1,..3} |x_i\, y_j - x_j\, y_i|^2. $$
Then, defining, for $i,j =1,..,3$, $i\neq j$,
$$ q_{ij}^f (v,w) := (v_i - w_i)\, \left(\frac{\pa_j f(v)}{f(v)} -  \frac{\pa_j f(w)}{f(w)} \right)
 - (v_j - w_j)\, \left(\frac{\pa_i f(v)}{f(v)} -  \frac{\pa_i f(w)}{f(w)} \right), $$
 we see that (remember (\ref{mded}))
$$  D_{x \mapsto |x|^{-1}}(f) = \frac12\int\int_{\R^3 \times \R^3} f(v)\,f(w)\,|v-w|^{-1}\, \left( \frac{\nabla f(v)}{f(v)} - \frac{\nabla f(w)}{f(w)} \right)^T\,$$
$$ \left(Id - \frac{(v-w)\otimes(v-w)}{|v-w|^2} \right)\, 
\left( \frac{\nabla f(v)}{f(v)} - \frac{\nabla f(w)}{f(w)} \right) \, dv dw $$
$$ = \frac14\, \sum_{i,j=1,..,3; i \neq j} \int\int_{\R^3 \times \R^3} f(v)\,f(w)\,|v-w|^{-3}
\, \left| q_{ij}^f (v,w) \right|^2 \, dv dw .$$
Expanding $q_{ij}^f$, we get (for $i,j =1,..,3$, $i\neq j$)
$$ q_{ij}^f (v,w)  = \bigg[v_i\, \frac{\pa_j f(v)}{f(v)} - v_j\, \frac{\pa_i f(v)}{f(v)} \bigg]
+ w_j\, \frac{\pa_i f(v)}{f(v)} - w_i\, \frac{\pa_j f(v)}{f(v)}$$  
$$ - v_i\, \frac{\pa_j f(w)}{f(w)} + v_j\, \frac{\pa_i f(w)}{f(w)} +
\bigg[w_i\, \frac{\pa_j f(w)}{f(w)} - w_j\, \frac{\pa_i f(w)}{f(w)} \bigg]. $$ 

 Choosing some $i,j=1,..,3$ with $i\neq j$, multiplying the above identity by $w_j\, (1+|w|^2)^{-3/4}\, f(w)$, and integrating w.r.t. the variable $w$ on $\R^3$,
 we end up with (remember that $f$ is assumed to be radially symmetric)
$$ \bigg( \int w_j^2\, (1+|w|^2)^{-3/4} \, f(w)\,dw \bigg)\, \frac{\partial_i f}{f}(v)$$
$$ =  v_i 
\bigg( \int  \bigg[(1+|w|^2)^{-3/4} - \frac32\, w_j^2\, (1+|w|^2)^{-7/4}\bigg] \, f(w)\,dw \bigg)\, $$
$$ + \int w_j\, (1+|w|^2)^{-3/4}  \,q_{ij}^f (v,w) \, f(w)\, dw . $$
As a consequence,
$$ \int \left| \frac{\partial_i f}{f}(v)\right|^2 \, f(v)\, (1+|v|^2)^{-3/2} dv $$
$$ \le 2  \bigg(\int v_i^2\,  (1+|v|^2)^{-3/2}\, f(v)\,dv\bigg)\,\, \bigg( \int w_j^2\, (1+|w|^2)^{-3/4} \, f(w)\,dw \bigg)^{-2}$$
$$
\times\, \bigg( \int  \bigg[(1+|w|^2)^{-3/4} - \frac32\, w_j^2\, (1+|w|^2)^{-7/4}\bigg] \, f(w)\,dw \bigg)^2 $$
$$+ \,2 \,\bigg( \int w_j^2\, (1+|w|^2)^{-3/4} \, f(w)\,dw \bigg)^{-2} $$
$$ \times\, \int f(v)\, (1+|v|^2)^{-3/2}\, \bigg( \int |w_j|\, (1+|w|^2)^{-3/4} \, |q_{ij}^f (v,w)|\, f(w)\, dw\bigg)^2\, dv $$
$$ \le \bigg( \int w_j^2\, (1+|w|^2)^{-3/4} \, f(w)\,dw \bigg)^{-2} \, \bigg(2 +  2 \int\int f(v)\,f(w)\, |v-w|^{-3} \, |q_{ij}^f(v,w)|^2\, dv dw $$
$$ \times  \sup_{v\in\R^3} \bigg[ (1+|v|^2)^{-3/2} \,\int f(w) \, |v-w|^3\, |w_j|^2\, (1+|w|^2)^{-3/2}\, dw \, \bigg]\, \bigg), $$
since (thanks to the radial symmetry) $\int v_i^2\,  (1+|v|^2)^{-3/2}\, f(v)\,dv \le 1$, and
$$ \bigg|\int  \bigg[(1+|w|^2)^{-3/4} - \frac32\, w_j^2\, (1+|w|^2)^{-7/4}\bigg] \, f(w)\,dw
\bigg|$$
$$ = \bigg|\int  \frac{1+ (\sum_{k\neq j} w_k^2) - \frac12\, w_j^2}{ (1+|w|^2)^{7/4} }\, f(w)\,dw \bigg| \le 1. $$
Then, estimating
$$\int f(w) \, |v-w|^3\, |w_j|^2\, (1+|w|^2)^{-3/2}\, dw $$
$$\le 4 \,|v|^3 \int f(w) \, |w_j|^2\, (1+|w|^2)^{-3/2}\, dw + 4  \int f(w) \, |w|^3\, |w_j|^2\, (1+|w|^2)^{-3/2}\, dw $$
$$\le \frac43\, |v|^3 + \frac{16}3, $$
we see that
$$ \int \left| \frac{\partial_i f}{f}(v)\right|^2 \, f(v)\, (1+|v|^2)^{-3/2} dv $$
$$ \le  \bigg( \int w_j^2\, (1+|w|^2)^{-3/4} \, f(w)\,dw \bigg)^{-2}
 \, \left(2 + \frac{128}3  \,D_{x \mapsto |x|^{-1}}(f) \right) . $$

We finally observe that (for all $0<\delta<R$, and $K>1$)
$$ \int w_j^2\, (1+|w|^2)^{-3/4} \, f(w)\,dw = \frac13\, \int |w|^2\,(1+|w|^2)^{-3/4} \, f(w)\,dw $$
$$ \ge  \frac13\, \int_{\delta \le |w| \le R} \delta^2\, (1 + R^2)^{-3/4} \, f(w)\, dw $$
$$ \ge \frac13\, \delta^2\, (1 + R^2)^{-3/4} \bigg(1 - \int_{\delta \ge |w|} f(w)\, 1_{f(w) \le K} \, dw $$
$$ -  \int_{\delta \ge |w|} f(w)\, \frac{|\ln f(w)|}{\ln K}  \, dw  - R^{-2} \int |w|^2  \, f(w)\,dw \bigg)$$
$$ \ge \frac13\, \delta^2\, (1 + R^2)^{-3/4} \, \bigg( 1 - \frac43\pi\delta^3\,K  - \frac{\bar{H}}{\ln K} - \frac3{R^2} \bigg). $$
Taking $R= \sqrt{12}$, $K =  \exp( 4\, \bar{H})$ and $\delta^3= \frac{3}{16\,\pi} \, \exp(-4\, \bar{H}) $, we end up with
$$ \int w_j^2\, (1+|w|^2)^{-3/4} \, f(w)\,dw \ge  \frac{(13)^{-3/4}}{12}\, \left(\frac{3}{16\pi}\right)^{2/3}\, \exp \left( - \frac{8}{3}\, \bar{H} \right), $$
which ends the proof of Theorem \ref{princ2} in the case when $f$ is radially symmetric and satisfies the
normalization (\ref{norma}).
$\square$

Note that an expression similar to $\exp \left( - \frac{8}{3}\, \bar{H} \right)$ appears 
in the proof of Theorem 1 of \cite{advw}, which can be seen as a local version of the same result, but for the Boltzmann equation.



\subsection{Proof of our main result in the general case} \label{sub2}

We begin by establishing two lemmas, which can be seen as a quantitative version of the
statement that it is impossible for $f$
to be concentrated on an hyperplane of $\R^N$, if it has a bounded entropy.

\begin{Lem} \label{lem1}
Let $f := f(v)\ge 0$ belong to $L^1_2(\R^N)$ and be
such that $\int_{\R^N} f(v)\, dv =1$, $\int_{\R^N} f(v)\, v\, dv =0$,
 and $\int_{\R^N} f(v)\, |v|^2\, dv =N$.
\smallskip

 Then, for all $i,j \in \N$ such that $i\neq j$, all $S>0$, and all $\lambda>0$,
$$  -{\hbox{ Det }} \left( \int_{\R^N} e^{- \lambda y^2}  \left[ \begin{array}{ccc}
 1  & y_j & y_i \\
 y_i  & y_i\,y_j & y_i^2 \\
 y_j & y_j^2& y_i\,y_j\end{array} \right]  \,f(y)\, dy \right) \, $$
 $$ \ge \frac12\, e^{-2N\,\lambda} \,\bigg( \int_{\R^N} e^{- \lambda y^2} \,f(y)\,y_i^2\,dy \,
\int_{\R^N} e^{- \lambda y^2} \,f(y)\,y_j^2\,dy - \bigg[  \int_{\R^N} e^{- \lambda y^2} \,f(y)\,y_i\,y_j\,dy \bigg]^2 \,\bigg)$$
$$ - \,4N \, (S\, (1 - e^{- \lambda\,S^2}) +2N\, S^{-1})^2. $$
 \end{Lem}

{\bf{Proof of Lemma \ref{lem1}}} : Expanding the determinant w.r.t. the first line, we get
(for any $i,j \in \N$ such that $i\neq j$ and all $\lambda>0$),
$$ - {\hbox{ Det }} \left( \int_{\R^N} e^{- \lambda y^2}  \left[ \begin{array}{ccc}
 1  & y_j & y_i \\
 y_i  & y_i\,y_j & y_i^2 \\
 y_j & y_j^2& y_i\,y_j\end{array} \right]  \,f(y)\, dy \right) \, $$
$$ = \bigg( \int_{\R^N} e^{- \lambda y^2} \,f(y)\,dy  \bigg) \, \bigg( \int_{\R^N} e^{- \lambda y^2} \,f(y)\,y_i^2\,dy \,
\int_{\R^N} e^{- \lambda y^2} \,f(y)\,y_j^2\,dy$$
$$  - \bigg[  \int_{\R^N} e^{- \lambda y^2} \,f(y)\,y_i\,y_j\,dy \bigg]^2 \,\bigg)$$
$$ + \bigg( \int_{\R^N} e^{- \lambda y^2} \,f(y)\,y_j\,dy  \bigg) \,  \bigg( \int_{\R^N} e^{- \lambda y^2} \,f(y)\,y_i\,dy \,
\int_{\R^N} e^{- \lambda y^2} \,f(y)\,y_i\,y_j\,dy$$
$$  -  \int_{\R^N} e^{- \lambda y^2} \,f(y)\,y_j\,dy \,
 \int_{\R^N} e^{- \lambda y^2} \,f(y)\,y_i^2\,dy \,\bigg)$$
$$ - \bigg( \int_{\R^N} e^{- \lambda y^2} \,f(y)\,y_i\,dy  \bigg) \,  \bigg( \int_{\R^N} e^{- \lambda y^2} \,f(y)\,y_i\,dy \,
\int_{\R^N} e^{- \lambda y^2} \,f(y)\,y_j^2\,dy$$
$$  -  \int_{\R^N} e^{- \lambda y^2} \,f(y)\,y_j\,dy \,
 \int_{\R^N} e^{- \lambda y^2} \,f(y)\,y_i\,y_j\,dy \,\bigg) $$
$$ \ge \bigg( \int_{\R^N} e^{- \lambda y^2} \,f(y)\,dy  \bigg) \, \bigg( \int_{\R^N} e^{- \lambda y^2} \,f(y)\,y_i^2\,dy \,
\int_{\R^N} e^{- \lambda y^2} \,f(y)\,y_j^2\,dy$$
$$  - \bigg[  \int_{\R^N} e^{- \lambda y^2} \,f(y)\,y_i\,y_j\,dy \bigg]^2 \,\bigg)$$
$$ - N\, \bigg( \bigg|\int_{\R^N} e^{- \lambda y^2} \,f(y)\,y_i\,dy\bigg| + \bigg|\int_{\R^N} e^{- \lambda y^2} \,f(y)\,y_j\,dy\bigg|\bigg)^2. $$
 Then, we observe (since $\int_{B(0,\sqrt{2N})^c} f(y)\,dy \le (2N)^{-1} \int_{B(0,\sqrt{2N})^c} y^2\, f(y)\,dy \le 1/2$) that
 $$ \int_{\R^N} e^{- \lambda y^2} \,f(y)\,dy \ge e^{-2N\,\lambda}\,\int_{B(0,\sqrt{2N})} f(y)\,dy \ge \frac12\, e^{-2N\,\lambda}, $$
 and (for any $S>0$)
 $$ \bigg|\int_{\R^N} e^{- \lambda y^2} \,f(y)\,y_i\,dy \bigg| = \bigg|\int_{\R^N} (e^{- \lambda y^2} - 1) \,f(y)\,y_i\,dy \bigg|$$
 $$ \le \int_{B(0,S)} S\, (1 - e^{-\lambda\,S^2})\, \,f(y)\,dy  + 2\,  \int_{B(0,S)^c}  |y_i| \,f(y)\,dy$$
 $$ \le S\, (1 - e^{-\lambda\,S^2}) + 2N\,S^{-1}. $$
 This ends the proof of Lemma \ref{lem1}.
 $\square$

\begin{Lem} \label{lem2} Let $f := f(v)\ge 0$ belonging to $L^1_2(\R^N) \cap L\,\ln L(\R^N)$ be
 such that $\int_{\R^N} f(v)\, dv =1$, $\int_{\R^N} f(v)\, v\, dv =0$,
 $\int_{\R^N} f(v)\, |v|^2\, dv =N$, and
 $\int_{\R^N} f(v)\, |\ln f(v)|\, dv \le \bar{H}$, where $\bar{H}>0$.  
\smallskip

 Then there exists a constant $\lambda_0 := \lambda_0(N,\bar{H})>0$ depending only on $N,\bar{H}$
 such that for all $i,j \in \N$, $i\neq j$  and $\lambda \in ]0, \lambda_0]$,
$$ \bigg| {\hbox{ Det }} \left( \int_{\R^N} e^{- \lambda y^2}  \left[ \begin{array}{ccc}
 1  & y_j & y_i \\
 y_i  & y_i\,y_j & y_i^2 \\
 y_j & y_j^2& y_i\,y_j\end{array} \right]  \,f(y)\, dy \right) \,\bigg| \ge  C(N, \bar{H}),$$
where $C := C(N, \bar{H})>0$ is also a constant which only depends on $N$ and $\bar{H}$.
 \end{Lem}

{\bf{Proof of Lemma \ref{lem2}}} : 
We observe that (for all $i,j \in \N$, $i\neq j$, and $\lambda >0$)
$$ \bigg(\int e^{-\lambda\,y^2} y_i^2 \,f(y)\,  dy \bigg)\, \bigg(\int e^{-\lambda\,y^2} y_j^2 \,f(y)\,  dy \bigg) 
 - \bigg(\int e^{-\lambda\,y^2} y_i\,y_j \,f(y)\,  dy \bigg)^2 $$
 $$ = \bigg(\int e^{-\lambda\,y^2} y_i^2 \,f(y)\,  dy \bigg)\, \bigg(\int e^{-\lambda\,y^2} \bigg|y_j 
- y_i \, \frac{\int e^{-\lambda\,z^2}\, z_iz_j\, f(z) dz}{\int e^{-\lambda\,z^2}\, z_i^2\, f(z) dz}\bigg|^2\, \,f(y)\,  dy \bigg)$$
$$ = \bigg(\int e^{-\lambda\,y^2} y_j^2 \,f(y)\,  dy \bigg)\, \bigg(\int e^{-\lambda\,y^2} \bigg|y_i 
- y_j \, \frac{\int e^{-\lambda\,z^2}\, z_iz_j\, f(z) dz}{\int e^{-\lambda\,z^2}\, z_j^2\, f(z) dz}\bigg|^2\, \,f(y)\,  dy \bigg).$$
Then, thanks to Cauchy-Schwarz (or Young) inequality, 
$$ \int e^{-\lambda\,z^2}\, z_iz_j\, f(z)\, dz \le \sup_{k=i,j} \int e^{-\lambda\,z^2}\, z_k^2\, f(z)\, dz, $$
so that
 $$ \bigg(\int e^{-\lambda\,y^2} y_i^2 \,f(y)\,  dy \bigg)\, \bigg(\int e^{-\lambda\,y^2} y_j^2 \,f(y)\,  dy \bigg) 
 - \bigg(\int e^{-\lambda\,y^2} y_i\,y_j \,f(y)\,  dy \bigg)^2 $$
 $$ \ge \inf_{|\theta |\le 1} \,\, \inf_{k=i,j;\, l=i,j;\,  k\neq l} \bigg(\int e^{-\lambda\,y^2} |y_k 
- \theta\, y_l |^2\, \,f(y)\,  dy \bigg)^2.$$ 
\smallskip

Using this estimate and Lemma~\ref{lem1}, we get for any $S>0$,
$$ \Gamma_{\lambda, i,j}(f) := - {\hbox{ Det }} \left( \int_{\R^N} e^{- \lambda y^2}  \left[ \begin{array}{ccc}
 1  & y_j & y_i \\
 y_i  & y_i\,y_j & y_i^2 \\
 y_j & y_j^2& y_i\,y_j\end{array} \right]  \,f(y)\, dy \right) \, $$
 $$ \ge \frac12\,e^{-2N\lambda} \, 
 \,\, \bigg( \inf_{|\theta |\le 1}\,\, \inf_{k=i,j;\, l=i,j;\, k\neq l}
 \int_{\R^N} e^{- \lambda y^2} \,f(y)\, |y_k - \theta\, y_l|^2\, dy \bigg)^2$$
$$ - \,4N\, (S\, (1 - e^{-\lambda\,S^2}) + 2N\,S^{-1})^2.$$
 
 Then, for any $\delta>0$ (and for any $|\theta |\le 1, \,\, k=i,j;\, l=i,j;\, k\neq l$),
 $$ \int_{\R^N} e^{- \lambda y^2} \,f(y)\, |y_k - \theta\, y_l|^2\, dy\ge \delta^2 \int_{|y_k - \theta\, y_l| \ge \delta}
 e^{- \lambda y^2} \,f(y)\, dy $$
 $$ \ge \delta^2\,  e^{- 2N\,\lambda} \int_{|y_k - \theta\, y_l| \ge \delta, \, |y|\le \sqrt{2N}} \,f(y)\, dy $$
 $$ = \delta^2\,  e^{- 2N\,\lambda} \bigg( 1 - \int_{|y_k - \theta\, y_l| \le \delta, \, |y|\le \sqrt{2N}} f(y)\, dy - 
\int_{ |y|\ge \sqrt{2N} } f(y)\, dy \bigg) $$
$$ \ge \delta^2\,  e^{- 2N \lambda} \bigg(   \frac{1}{2} - \int_{|y_k - \theta\, y_l| \le \delta, \, |y|\le \sqrt{2N} } f(y)\, dy
  \bigg). $$
\smallskip

 As a consequence, for any $K>1$, $S>0$, $\delta>0$,
 
 $$ \Gamma_{\lambda,i,j}(f) \ge \frac12\,e^{-6N\lambda} \, \delta^4  \,  \bigg(   \frac{1}{2} - \sup_{|\theta |\le 1} \,\,\sup_{k=i,j;\, l=i,j;\, k\neq l}
\int_{|\frac{y_k - \theta\, y_l}{\sqrt{1+\theta^2}}| \le \frac{\delta}{\sqrt{1+\theta^2}}, \, |y|\le \sqrt{2N} } f(y)\, dy
  \bigg)^2 $$
$$- \, 4N\, (S\, (1 - e^{-\lambda\,S^2}) + 2N\,S^{-1})^2 $$
$$ \ge \frac12\,e^{-6N\lambda} \, \delta^4  \,  \bigg(\frac12 - 2K\,(2 \sqrt{2N})^{N-1}\,\delta - \frac{\bar{H}}{\ln K} \bigg)^2
-\, 4 N\, (S\, (1 - e^{-\lambda\,S^2}) + 2N\,S^{-1})^2 .$$
Taking $K= e^{4\bar{H}}$, $\delta = 16^{-1}\,(8N)^{- \frac{N-1}2} \, e^{-4\bar{H}}$, we see that (for any $S>0$)
$$ \Gamma_{\lambda,i,j}(f) \ge \frac1{128} (16\,(8N)^{\frac{N-1}2})^{-4} \, e^{-16 \bar{H}}\, e^{-6N\lambda} 
- \, 4N\, (2S^2\, (1 - e^{-\lambda\,S^2})^2 + 8N^2\,S^{-2}) . $$
Finally, taking $S>0$ in such a way that 
$$ 32 N^3\,S^{-2} = 2^{-24 - 6(N-1)}\,N^{-2(N-1)}\, e^{-16 \bar{H}}, $$
we get
$$ \Gamma_{\lambda,i,j}(f) \ge 2^{-17 - 6N}\, N^{-2(N-1)}\, e^{-16 \bar{H}} \, \left(e^{-6N\lambda} - \frac12 \right) $$
 $$ - 2^{26 + 6N} \, N^{2(N+1)} \, e^{16 \bar{H}} \, (1 - e^{- 2^{23+6N} \, N^{2N+1} \, e^{16 \bar{H}} \,\lambda})^2 .$$
The limit of this formula when $\lambda \to 0$ is
$2^{-18 - 6N}\, N^{-2(N-1)}\, e^{-16 \bar{H}}$.
This ends the proof of Lemma \ref{lem2}.
\smallskip

Note that when $N=3$, Lemma \ref{lem2} can be made numerically explicit:  $\Gamma_{\lambda,i,j}(f) \ge 2^{-38}\, 3^{-4} \, \exp ( -16 \bar{H} )$,
provided that $0 < \lambda \le \lambda_0 = 2^{-82}\, 3^{-13}  \, \exp ( -24 \bar{H} )$. 
$\square$
\medskip

We now turn to the 
\medskip

{\bf{Proof of Theorem \ref{princ1}}}: 
We shall assume in this proof, without loss of generality, that 
\begin{equation} \label{norm}
 \int_{\R^N} f(v)\, dv = 1, \quad \int_{\R^N} f(v)\, v\, dv =0, 
\quad \int_{\R^N} f(v)\, |v|^2
\, dv =N.
\end{equation}
It indeed amounts to a simple change of unknown and variables of the form
\begin{equation} \label{lescaling}
f(v) \longrightarrow a \, f(b\,v+c), \qquad a,b\in\R, c\in\R^N.
\end{equation}
In the computations below, we do not take care of
the points where $f =0$. We explain how to justify 
the computations where $\frac{1}{f}$ appears at the end of the proof of Theorem \ref{princ1}.
\smallskip

We first observe that (for all $x,y\in\R^N$)
$$ y^T\,(|x|^2 \, Id- x \otimes x)\,y = \frac12 \, \sum_{i,j=1,..N} |x_i\, y_j - x_j\, y_i|^2. $$
Then, defining, for $i,j =1,..,N$, $i\neq j$,
$$ q_{ij}^f (v,w) := (v_i - w_i)\, \left(\frac{\pa_j f(v)}{f(v)} -  \frac{\pa_j f(w)}{f(w)} \right)
 - (v_j - w_j)\, \left(\frac{\pa_i f(v)}{f(v)} -  \frac{\pa_i f(w)}{f(w)} \right), $$
 we see that (remember (\ref{mded}))
$$  D_{\psi}(f) = \frac12\int\int_{\R^N \times \R^N} f(v)\,f(w)\,\psi(|v-w|)\, \left( \frac{\nabla f(v)}{f(v)} - \frac{\nabla f(w)}{f(w)} \right)^T\,$$
$$ \left(Id - \frac{(v-w)\otimes(v-w)}{|v-w|^2} \right)\, 
\left( \frac{\nabla f(v)}{f(v)} - \frac{\nabla f(w)}{f(w)} \right) \, dv dw $$
$$ = \frac14\, \sum_{i,j=1,..N} \int\int_{\R^N \times \R^N} f(v)\,f(w)\,\frac{\psi(|v-w|)}{|v-w|^2}\,
\left| q_{ij}^f (v,w) \right|^2 \, dv dw .$$
Expanding $q_{ij}^f$, we get
$$ q_{ij}^f (v,w)  = \bigg[v_i\, \frac{\pa_j f(v)}{f(v)} - v_j\, \frac{\pa_i f(v)}{f(v)} \bigg]
+ w_j\, \frac{\pa_i f(v)}{f(v)} - w_i\, \frac{\pa_j f(v)}{f(v)}$$  
$$ - v_i\, \frac{\pa_j f(w)}{f(w)} + v_j\, \frac{\pa_i f(w)}{f(w)} +
\bigg[w_i\, \frac{\pa_j f(w)}{f(w)} - w_j\, \frac{\pa_i f(w)}{f(w)} \bigg]. $$ 
We now select some $i,j =1,..,N$, $i\neq j$, and integrate the identity above
against $e^{- \lambda w^2}\, f(w)\,dw$, for some $\lambda>0$. We get
$$ \int_{\R^N} q_{ij}^f (v,w)\, e^{- \lambda w^2}\, f(w)\,dw  = \bigg[v_i\, \frac{\pa_j f(v)}{f(v)} - v_j\, \frac{\pa_i f(v)}{f(v)} \bigg] \, \bigg( \int_{\R^N} e^{- \lambda w^2}\, f(w)\,dw  \bigg)$$
$$ + \bigg( \int_{\R^N} w_j\, e^{- \lambda w^2}\, f(w)\,dw \bigg)\,\frac{\pa_i f(v)}{f(v)}
 - \bigg(\int_{\R^N} w_i\, e^{- \lambda w^2}\, f(w)\,dw  \bigg)\,\frac{\pa_j f(v)}{f(v)}$$  
$$ - 2\lambda \,v_i\, \bigg( \int_{\R^N} w_j \, e^{- \lambda w^2}\, f(w)\,dw \bigg) 
 + 2\lambda \,v_j\, \bigg( \int_{\R^N} w_i \, e^{- \lambda w^2}\, f(w)\,dw \bigg) . $$ 
We then integrate it against $w_k \,e^{- \lambda w^2}\, f(w)\,dw$, for some $k=1,..,N$ (and some $\lambda>0$)
and get
$$ \int_{\R^N} q_{ij}^f (v,w)\, w_k\,e^{- \lambda w^2}\, f(w)\,dw  = \bigg[v_i\, \frac{\pa_j f(v)}{f(v)} - v_j\, \frac{\pa_i f(v)}{f(v)} \bigg] \, \bigg( \int_{\R^N} w_k \,e^{- \lambda w^2}\, f(w)\,dw \bigg) $$
$$ + \bigg( \int_{\R^N} w_j\, w_k\, e^{- \lambda w^2}\, f(w)\,dw \bigg) \, \frac{\pa_i f(v)}{f(v)}
 - \bigg( \int_{\R^N} w_i\, w_k\, e^{- \lambda w^2}\, f(w)\,dw \bigg)\, \frac{\pa_j f(v)}{f(v)}  $$
$$+ \,v_i\, \bigg( \int_{\R^N} (\delta_{jk} -2\lambda\,w_j\,w_k) \, e^{- \lambda w^2}\, f(w)\,dw
\bigg) 
 - v_j\, \bigg( \int_{\R^N} (\delta_{ik} -2\lambda\,w_i\,w_k) \, e^{- \lambda w^2}\, f(w)\,dw \bigg)   $$
$$ +\,   \int_{\R^N} (\delta_{ik}\, w_j  - \delta_{jk}\,w_i) \, e^{- \lambda w^2}\, f(w)\,dw  . $$
Using $k=i,j$,
 and considering the above identities as a $3 \times 3$ system for the unknowns
$v_i\, \frac{\pa_j f(v)}{f(v)} - v_j\, \frac{\pa_i f(v)}{f(v)}$, $\frac{\pa_i f(v)}{f(v)}$ and
$\frac{\pa_j f(v)}{f(v)}$, Cramer formulas enable the computation of $\frac{\pa_i f(v)}{f(v)}$ in terms of $q^f_{ij}$:
$$ \frac{\pa_i f(v)}{f(v)} = \frac{ Det \, \bigg(  \int_{\R^N} e^{- \lambda w^2}\, f(w)
 \left[ \begin{array}{ccc}
 1 & Z_1(f)(v,w) & w_i \\
 w_i & Z_2(f)(v,w) &  w_i^2 \\
 w_j & Z_3(f)(v,w) &  w_i\,w_j \end{array}\, \right] \,dw\, \bigg)} { Det \, \bigg(  \int_{\R^N} e^{- \lambda w^2}\, f(w)
 \left[ \begin{array}{ccc}
  1 & w_j &  w_i \\
 w_i & w_j\,w_i &  w_i^2 \\
 w_j & w_j^2 &  w_i\,w_j \end{array} \, \right] \,dw\, \bigg)}\, , $$
 where
$$  Z_1(f)(v,w) =  q_{ij}^f (v,w)  + 2\lambda \,v_i\, w_j  - 2\lambda v_j\, w_i, $$
 $$ Z_2(f)(v,w) =  q_{ij}^f (v,w)\, w_i  - v_i\, ( -2\lambda\,w_j\,w_i)  + v_j\, (1 -2\lambda\,w_i^2) -  w_j, $$
$$ Z_3(f)(v,w) = q_{ij}^f (v,w)\, w_j - v_i\, (1 -2\lambda\,w_j^2) + v_j\, ( -2\lambda\,w_i\,w_j) +  w_i . $$
From now on, we denote by $C_1(..), C_2(..), ...$ various constants depending only on the parameters
indicated in the parenthesis.
\medskip

 Using Lemma \ref{lem2} (and the notation $\Gamma_{\lambda,i,j}$ introduced in its proof), we see that, taking $\lambda = \lambda_0(N, \bar{H})$ (and denoting it by $\lambda_0$),
 $$ \bigg| \frac{\pa_i f(v)}{f(v)} \bigg| \le (\Gamma_{\lambda_0}(f))^{-1}\, 
 2N^{3/2} \,  \int_{\R^N} e^{- \lambda_0\, w^2}\, f(w)$$
$$ \times\, \bigg( |Z_1(f)(v,w)| +  |Z_2(f)(v,w)| +  |Z_3(f)(v,w)| \bigg)\, dw $$
 $$ \le C(N, \bar{H})^{-1}\, C_1(N,\bar{H})\, \bigg( [1+ |v|] + \int_{\R^N} q_{ij}^f (v,w)\, (1+|w|^2)^{1/2}\, e^{- \lambda_0 w^2}\, f(w)\,dw\bigg).$$
 Then
 $$ \int_{\R^N} f(v)\, \bigg| \frac{\pa_i f(v)}{f(v)} \bigg|^2\, (1+|v|^2)^{\inf(\gamma_1/2, -1)}\, dv $$
 $$ \le C_2(N, \bar{H})\, \bigg( \int_{\R^N} f(v)\, (1+|v|^2)^{\inf(1+\gamma_1/2, 0)} \, dv $$
$$+  \int_{\R^N} f(v)\, (1+|v|^2)^{\inf(\gamma_1/2, -1)}\, \bigg|\int_{\R^N} q_{ij}^f (v,w)\,(1+|w|^2)^{1/2}\, e^{- \lambda_0 w^2}\, f(w)\,dw\bigg|^2 \, dv\bigg) $$
$$ \le C_2(N, \bar{H})\, \bigg(1 + \int_{\R^N} \int_{\R^N} f(v)\, f(w)\, \frac{\psi(|v-w|)}{|v-w|^2} \left|q_{ij}^f (v,w)\right|^2\,
dvdw $$
$$ \times\,  \sup_{v\in\R^N} \bigg[  (1+|v|^2)^{\inf(\gamma_1/2,-1)}\, \int_{\R^N} \frac{|v-w|^2}{\psi(|v-w|)}\,  (1+|w|^2)\, e^{- 2\lambda_0 w^2}\,f(w)\, dw  \bigg]\,\bigg) $$
$$ \le C_3(N, \bar{H}, K_3)\, \bigg(1 + D_{\psi}(f)\, \sup_{v\in\R^N}   \bigg[\,(1+|v|^2)^{\inf(\gamma_1/2, -1)}$$
$$ \times\,  \int_{\R^N} \sup \bigg(|v-w|^{-\gamma_1}, |v-w|^2 \bigg)\,  (1+|w|^2)\, e^{- 2\lambda_0 w^2}\,f(w)\, dw \,\bigg]\, \bigg)$$
$$ \le C_4(N, \bar{H}, K_3, \gamma_1) \, (1 + D_{\psi}(f)). $$ 
This concludes the proof of Theorem \ref{princ1}, up to the treatment of the points where $f=0$.
\medskip

We briefly explain here how to treat this problem. First, one needs to start from the definition of 
$D_{\psi}(f)$ appearing in the last line of (\ref{mded}), which requires no extra assumption on $f$ (other than
$f \in L^1_{loc}$). Then, one uses $g= \sqrt{f}$, and the integration
against $e^{- \lambda\,w^2} \, f(w)\, dw$, $e^{- \lambda\,w^2} \, w_i\, f(w)\, dw$ is replaced
by an integration
against $e^{- \lambda\,w^2} \, g(w)\, dw$, $e^{- \lambda\,w^2} \, w_i\, g(w)\, dw$. Though such a presentation 
is more rigorous, we did not adopt it
 in our proof in order to keep working with $f$ (and not
$g = \sqrt{f}$) and the
definition of $D_{\psi}(f)$ appearing in the second line of (\ref{mded}) (and not the one appearing in the last line), since those quantities are much 
more familiar to specialists of kinetic theory. 
$\square$ 
 \medskip
 
 We now turn to the 
\medskip

{\bf{Proof of Theorem \ref{princ2}}}:  The first part is a direct consequence of Theorem 
\ref{princ1}, when $N=3$, $\gamma_1 = -3$, $K_3=1$. The second part (case of radially symmetric $f$) was already proven in subsection \ref{sub1}.
 $\square$
 
 \section{$H$-solutions of the Landau equation are weak solutions} \label{weaksose}
 
 Though a large part of the material of this section is either classical or analogous to
 well-known results on the Boltzmann equation (cf. the remarks at the end of the section),
 we provide detailed proofs for the sake of completeness.

 \subsection{Weighted $L^p$ estimates}

We begin this section with a direct application of Sobolev inequalities in the whole space. We write a general lemma (the dimension $N$ as well as the exponent of the weight are arbitrary)
linking the weighted $H^1$ norm of $\sqrt{f}$ to a weighted $L^{\frac{N}{N-2}}$ norm of $f$.

\begin{Lem} \label{lemm1} 
Let $N \in \N$, $N\ge 3$, and $\gamma_1 \in \R$.
Then, there exists a constant depending only on $N$ and $\gamma_1$ (denoted by $C:= C(N, \gamma_1)>0$)
 such that for all $f:= f(v)\ge 0$ lying in $L^1_2(\R^N)$
$$  \bigg( \int_{\R^N} |f(v)|^{\frac{N}{N-2}}  \, (1 + |v|^2)^{\frac{N}{N-2}\, \inf(\frac{\gamma_1}{2}, -1) }\, dv \bigg)^{\frac{N-2}{N}} $$
$$ \le C \,\bigg[ \int_{\R^N} f(v)\, (1+|v|^2)\, dv  +   \int_{\R^N} |\nabla\sqrt{f(v)}|^2 \, 
(1 + |v|^2)^{\inf(\frac{\gamma_1}{2}, -1)}\, dv  \bigg]. $$
In the inequality above, the right-hand side is considered as equal to $+\infty$ if $\sqrt{f}$
is not in $H^1_{loc}(\R^N)$.
\medskip
 
 In particular, for $N=3$ and $\gamma_1= -3$ (that is, the Coulomb case),
 $$ \bigg(\int_{\R^3} |f(v)|^3 \, (1 + |v|^2)^{- \frac92} \, dv\bigg)^{\frac{1}3} $$
$$ \le \frac6{\sqrt{\pi}} \int_{\R^3} f(v)\,dv
 +   \frac8{3\,\sqrt{\pi}} \int_{\R^3} |\nabla\sqrt{f(v)}|^2 \, (1 + |v|^2)^{\frac32} dv. $$
 \medskip
 
 If $N= 2$, and $\gamma_1 \in \R$, then for any $q\ge 1$,  there exists a constant depending only on $q$ and $\gamma_1$ (denoted by $C:= C(q, \gamma_1)>0$)  such that 
 for all $f\ge 0$ lying in $L^1_2(\R^2)$,
 $$  \bigg(\int_{\R^2} |f(v)|^{q}  \, (1 + |v|^2)^{q\, \inf(\frac{\gamma_1}2, -1)}\, dv \bigg)^{\frac1q}$$
$$ \le C \,\bigg[  \int_{\R^2} f(v)\,
 dv  + \int_{\R^2} |\nabla\sqrt{f(v)}|^2 \, 
(1 + |v|^2)^{\inf(\frac{\gamma_1}2, -1)}\, dv \bigg]. $$
 \end{Lem}

{\bf{Proof of Lemma \ref{lemm1}}} : We denote by $C(N)$ (resp. $C(N, \gamma_1)$) any constant depending only on $N$ (resp. $N$, $\gamma_1$). 
We recall the standard Sobolev inequality for functions of $H^1(\R^N)$
 (for $N\ge 3$):
 $$ \int_{\R^N} |h|^{\frac{2N}{N-2}}  \le C(N)\, \bigg( \int_{\R^N} |\nabla h|^2 \bigg)^{\frac{N}{N-2}} . $$
Applying it to $h(v) = g(v)\, (1+|v|^2)^{\frac12\, \inf(\gamma_1/2, -1)}$, we see that
$$  \int_{\R^N} |g(v)|^{\frac{2N}{N-2}} \, (1 + |v|^2)^{\frac{N}{N-2} \, \inf(\frac{\gamma_1}{2}, -1)} \, dv 
\le C(N)\, \bigg( \int_{\R^N} |\nabla g(v)|^2 \, (1 + |v|^2)^{\inf(\frac{\gamma_1}{2},-1)} dv \bigg)^{\frac{N}{N-2}} $$
$$
+ \,C(N, \gamma_1)\, \bigg(  \int_{\R^N} |g(v)|^2
 \, dv \bigg)^{\frac{N}{N-2}} . $$
Then, taking $f = g^2$, we end up with
$$ \int_{\R^N} |f(v)|^{\frac{N}{N-2}}  \, (1 + |v|^2)^{\frac{N}{N-2}\, \inf(\frac{\gamma_1}{2}, -1)}\, dv  
\le C(N) \,\bigg(  \int_{\R^N} |\nabla\sqrt{f(v)}|^2 \, (1 + |v|^2)^{\inf(\frac{\gamma_1}2, -1)}\, dv  \bigg)^{\frac{N}{N-2}}$$
$$ + \,  C(N, \gamma_1) \,\bigg(  \int_{\R^N} f(v) \,
 dv  \bigg)^{\frac{N}{N-2}}. $$
\medskip

The Coulomb case is treated by following the constants in the above computations, and by using 
the best constant in the Sobolev embedding.
\medskip

In the case when $N=2$, we proceed as in the case when $N\ge 3$, but starting from 
the Sobolev inequality (which holds for all $q\in [1,+\infty[$, and with a constant $C(q)$ which
depends only on $q$): 
$$ \int_{\R^2} |h|^{2q}  \le C(q)\, \bigg( \int_{\R^2} |\nabla h|^2 \bigg)^{q} . $$
$\square$



\subsection{Weak solutions} 

We also start with a lemma which holds for all dimensions, and any function $\psi$ bounded above in the vicinity of $0$ by a power corresponding to very soft potentials. It is a simple consequence of Young's inequality for convolutions.
 
\begin{Lem} \label{lemm2} 
Let $N\in \N$, $N\ge 1$.
We assume that $\psi$ is a function satisfying
$$ \forall z\ge 0, \qquad 0 \le \psi(|z|) \le  K_1\, |z|^2 + K_2\, |z|^{\gamma_2 + 2}, $$
for some $K_1,K_2>0$ and $\gamma_2 \in ]-4, -2[$. 
\smallskip

Then for any $R>0$, and $r < N/(-\gamma_2 -2)$, there exists a constant
$C:= C(r, N, K_1,K_2,\gamma_2)>0$ depending only on the parameters of $\psi$, $N$
and $r$, such that for any $f \in L^1_2(\R^N) \cap L^{r'}_{loc}(\R^N)$
 (with $\frac1r + \frac1{r'} =1$),
$$  \int_{\R^N}\int_{B(0,R)}  f(v)\, f(w)\, \psi(|v-w|) \, dv dw $$
$$  \le C\, ||f||_{L^1(\R^N)}  \, \left( ||f||_{L^1_2(\R^N)} 
+ || f\, 1_{|\,\cdot\,| \le R} ||_{L^{r'}(\R^N)} \right). $$
 In particular, if $N\ge 3$,  there exists a constant
$C:= C(N, K_1,K_2,\gamma_2)>0$ depending only on the parameters of $\psi$ and $N$, 
such that for any $f \in L^1_2(\R^N) \cap L^{\frac{N}{N-2}}_{loc}(\R^N)$,
 $$  \int_{\R^N}\int_{B(0,R)}  f(v)\, f(w)\, \psi(|v-w|) \, dv dw $$
$$  \le C  \, ||f||_{L^1(\R^N)}  \, \left(||f||_{L^1_2(\R^N)}
 + || f\, 1_{|\,\cdot\,| \le R} ||_{L^{\frac{N}{N-2}}(\R^N)} \right). $$
 \end{Lem}

{\bf{Proof of Lemma \ref{lemm2}}} :
$$  \int_{\R^N}\int_{B(0,R)}  f(v)\, f(w)\, \psi(|v-w|) \, dv dw $$
$$ \le \int_{\R^N}\int_{B(0,R)}  f(v)\, f(w)\, (K_1 + K_2\, |v-w|^{\gamma_2+2}) \, 1_{|v-w|\le 1}\, dv dw $$
$$ +  \int_{\R^N}\int_{B(0,R)}  f(v)\, f(w)\, (K_1\, |v-w|^{2} + K_2) \, 1_{|v-w|\ge 1}\, dv dw $$
$$ \le K_1\, ||f||_{L^1(\R^N)}^2  + K_2 \int_{\R^{N}}\int_{B(0,R)} f(v) \, f(w) \,\,
 |v-w|^{\gamma_2+2} \, 1_{|v-w|\le 1}\, \,dv dw $$
 $$ +\, 4\,K_1 \,||f||_{L^1(\R^N)}\,||f||_{L^1_2(\R^N)}  + K_2\, ||f||_{L^1(\R^N)}^2  $$
$$ \le \left(4\,K_1\,||f||_{L^1_2(\R^N)} + (K_1 + K_2)\, ||f||_{L^1(\R^N)} \right) \, ||f||_{L^1(\R^N)} $$
$$ +\,K_2\, || f\, 1_{|\,\cdot\,| \le R} ||_{L^{r'}(\R^N)} \, || f ||_{L^1(\R^N)} \,
 || x \mapsto |x|^{\gamma_2 +2}\,  1_{|\,\cdot\,| \le 1} ||_{L^r(\R^N)}, $$
 thanks to Young's inequality for convolutions. 
\smallskip

Note that the last term is finite when $r < N/(- \gamma_2 -2)$, which is equivalent to $r'> N/(N+ \gamma_2+2)$. Observing that (when $N\ge 3$)
 $N/(N+\gamma_2+2) < N/(N-2)$ since $\gamma_2>-4$, we see that one can take $r'= N/(N-2)$, 
 so that the particular case (at the end of the statement of Lemma \ref{lemm2}) also holds.
$\square$
\medskip

 We are now in a position to write down the
 \medskip
 
 {\bf{Proof of Corollary \ref{weakso}}}: We know (cf. \cite{vill:new:97}) that $H$-solutions $f: = f(t,v)$ of the Landau equation (in the Coulomb case) satisfy the {\it{a priori}} estimate (\ref{meeap}).
As a consequence, they lie in $L^{\infty}([0,T]; L^1_2(\R^3))$ and 
satisfy the boundedness of the entropy dissipation, that is  
 $$ \int_0^T D_{x \mapsto |x|^{-1}}(f(t,\cdot))\, dt < +\infty . $$
\smallskip

 Using the main result of this paper (that is, Theorem \ref{princ2}), we see that
 $$ \int_0^T \int_{\R^3} |\nabla \sqrt{f(t,v)}|^2\, (1+ |v|^2)^{- 3/2} \, dv \, dt < + \infty . $$
 Then, using Lemma \ref{lemm1}, we end up with
 $$ \int_0^T \bigg(\int_{\R^3} |f(t,v)|^3\, (1+ |v|^2)^{- 9/2} \, dv  \bigg)^{1/3}\, dt < + \infty , $$
 so that $f \in L^1([0,T]; L^3_{-3}(\R^3))$, which is the first statement of Corollary \ref{weakso}
\medskip 
 
 Then, we consider $\varphi \in C^2_c([0,T[ \times \R^3)$, and $R>0$ such that 
$$ \cup_{t\in [0,T[} {\hbox{ Supp }} f(t,\cdot) \subset B(0,R). $$
 Then (since $v$ and $w$ play the same role), for any $i,j =1,..,3$,
 $$\int_0^T \iint_{\R^3\times\R^3}
 f(t,v)\, f(t,w)\, a_{ij}(v-w)\, 
\Bigl |\partial_{ij} \varphi(t,v) + \partial_{ij}\varphi(t,w) \Bigr | \, dv dw\, dt $$
$$ \le 2 \, ||\partial_{ij} \varphi||_{L^{\infty}([0,T] \times \R^3)}\, \int_0^T \iint_{\R^3\times B(0,R)} f(t,v)\, f(t,w)\, |v-w|^{-1} \, dv dw\, dt $$
$$ \le C(\varphi)\,||f||_{L^{\infty}([0,T]; L^1(\R^3))}\, 
\left( T\,||f||_{L^{\infty}([0,T]; L^1_2(\R^3))} +  ||(t,v) \mapsto f(t,v)\, 1_{|v|\le R}||_{L^1([0,T]; L^3(\R^3))}
\right) , $$
thanks to Lemma \ref{lemm2} (particular case, with $N=3$), and with $C(\varphi)$ depending only on $||\nabla^2\varphi||_{L^{\infty}(\R^3)}$. 
This last quantity is finite since $f \in L^1([0,T]; L^3_{-3}(\R^3)$.
\smallskip

We recall that in the Coulomb case, for $i=1,..,3$, $|b_i(z)| \le 2\,|z|^{-2}$. 
Then (since $v$ and $w$ play the same role), taking $\varphi$ as above, and $i=1,..,3$,
$$  \int_0^T\iint_{\R^3\times\R^3} f(t,v)\, f(t,w)\, |b_i(v-w)|\, 
\Bigl| (\partial_i \varphi)(t,v) - (\partial_i \varphi)(t,w) \Bigr|\, dv dw\, dt$$
$$ \le 4 \,\int_0^T\iint_{\R^3\times B(0,R)} f(t,v)\, f(t,w)\, |v-w|^{-2}\, 
\Bigl| (\partial_i \varphi)(t,v) - (\partial_i \varphi)(t,w) \Bigr|\, dv dw\, dt$$
$$ \le 4 \,  \sup_{j=1,..,3} ||\partial_{ij} \varphi||_{L^{\infty}([0,T] \times \R^3)}\,
 \int_0^T\iint_{\R^3\times B(0,R)} f(t,v)\, f(t,w)\, |v-w|^{-1}\,  dv dw\, dt,$$
 which, once again, is finite thanks to Lemma \ref{lemm2} and the estimate $f \in L^{\infty}([0,T]; L^1_2(\R^3)) \cap L^1([0,T]; L^3_{-3}(\R^3))$.
 \medskip
 
 The last statement of Corollary \ref{weakso} is thus proven. 
$\square$
\medskip

{\bf{Remark}}: As stated earlier, the
 computations of this section  are either classical or direct extensions and variants of
classical computations. They are close, for example, to computations of \cite{advw} (section 7, second application) or \cite{fournier_guerin} (proof of Proposition 3.3, step 1).
Note that in the proof, we only need the estimate $f\, 1_{|v|\le R} \in L^1([0,T]; L^3(\R^3))$,
 and not $f\in L^1([0,T]; L^3_{-3}(\R^3))$, so that the ``local in $v$'' approach
of \cite{advw} is well adapted to get this result.
\medskip

Observing the assumptions of Theorem \ref{princ2} and Lemma \ref{lemm2},
we see that if 
$$ \forall z\ge 0, \qquad  K_3\, \inf(1, |z|^{\gamma_1 + 2})
\le \psi(|z|) \le  K_1\, |z|^2 + K_2\, |z|^{\gamma_2 + 2}, $$
for some $K_1,K_2, K_3>0$ and $\gamma_1<0$, $\gamma_2 \in ]-4, -2[$, then reasonable 
solutions $f$ of the Landau equation with initial data having a finite mass, energy and entropy,
naturally satisfy (for all $T>0$)
the estimate $\int_0^T\int_{\R^N} |\nabla \sqrt{f(t,v)}|^2\, (1+ |v|^2)^{\inf(\gamma_1/2, -1)} \, dv dt < + \infty$, and thus belong to  $L^1([0,T]; L^{\frac{N}{N-2}}_{\inf(\gamma_1, -2)}(\R^N))$ when $N\ge3$
(and $L^1([0,T]; L^{q}_{\inf(\gamma_1, -2)}(\R^2))$ for all $q \in [1,+\infty[$ when $N=2$). This is
sufficient to define weak solutions according to Lemma \ref{lemm2}, and the proof of Corollary~\ref{weakso}.
\medskip

We also wish to emphasize the appearance of the limiting case $\gamma =-4$ in Lemma~\ref{lemm2}. 
This was already observed in \cite{advw}, and, from our point of view, justifies the 
terminology of ``very soft potentials'' for the case when $\gamma \in ]-4,-2[$, rather than
$\gamma \in ]-3,-2[$ or $\gamma \in [-3,-2[$. Note also that the value $\gamma = -4$ suprisingly appears
for all dimensions $N$.

 \section{Estimates for moments} \label{momse}

We present here an estimate of propagation of moments which holds for any weak solutions of the Landau equation,
in arbitrary dimension, when the function $\psi$ is bounded above around $0$ by a power law of either very soft or (moderately) soft type. We also impose that $\psi$ is bounded below at infinity by an arbitrary power law, and that $\psi$ is bounded above at infinity by a (moderately) soft potential.
\par
This estimate implies Proposition \ref{momco} in the Coulomb case.
\medskip

Our result of propagation of moments writes
\medskip

 \begin{Prop}\label{mom1}
 Let $\psi$ satisfy
$$ \forall z\ge 0, \qquad K_3\, \inf(1, |z|^{\gamma_1 + 2}) \le \psi(z) \le  K_1\, |z|^{2-\delta} + K_2\, |z|^{\gamma_2+2}, $$
for some $K_1,K_2, K_3>0$ and $\delta \in ]0,2]$, $\gamma_1 \le 0$, $\gamma_2 \in ]-4, -\delta]$. 
\par
 Let $N\in \N$, $N\ge 2$, $T>0$, and $f:=f(t,v)\ge 0$ be a weak solution of the Landau equation (\ref{Landau}) on $[0,T]\times\R^N$
associated to $\psi$ and an initial datum $f_{in} \in L^1_2 \cap L\,\ln L (\R^N)$,
in the following sense:
We assume that $f \in C( [0,T]; {\mathcal{D}}'(\R^N)) \cap L^{\infty}([0,T]; L^1_2 \cap L\,\ln L (\R^N))$, satisfies the conservation of mass, momentum, energy,  
$$ \forall t \in [0,T], \qquad 
 \int_{\R^N} f(t,v) \, \left(\begin{array}{c} 1\\ v_i \\ |v|^2/2 \end{array} \right)\, dv
=  \int_{\R^N} f_{in}(v) \, \left(\begin{array}{c} 1\\ v_i \\ |v|^2/2 \end{array} \right)\, dv,$$
the (uniform in time) boundedness of the entropy 
$$ \sup_{t\in [0,T]} \int_{\R^N} f(t,v) \, |\ln f(t,v)| \, dv <+\infty, $$
and has a finite (integrated in time) $L^{\frac{N}{N-2}}_{ \inf(\gamma_1, -2)}$ norm, that is 
\begin{equation} \label{cutie}
 {\mathcal{Q}}_{T,N,\gamma_1}(f) := \int_0^T  \bigg( \int_{\R^N} |f(v)|^{\frac{N}{N-2}}  \,
 (1 + |v|^2)^{\frac{N}{N-2}\, \inf(\frac{\gamma_1}{2}, -1) }\, dv \bigg)^{\frac{N-2}N} \, dt < +\infty,
 \end{equation}
if $N\ge 3$ (in this assumption, $N/(N-2)$ is replaced by all $q\in [1,+\infty[$ if $N=2$).
\par
 Moreover, we assume that $f(0, \cdot) = f_{in}$ and $f$
 satisfies the (strong w.r.t. $t$) weak form of the Landau equation, that is,
for all $\displaystyle \varphi := \varphi(v) \in
C^2_c({\Bbb R}^N)$, and $t_1 \le t_2 \in [0,T]$,
\begin{equation} \label{starpri}
\int_{\R^N}  f(t_2,v)\, \varphi(v)\, dv - \int_{\R^N} f(t_1, v)\, \varphi(v)\, dv 
\end{equation}
$$ = \frac{1}{2} \sum_{i=1}^3  \sum_{j=1}^3 \int_{t_1}^{t_2} \iint_{\R^N\times\R^N}
 f(s,v)\, f(s,w)\, a_{ij}(v-w)\, 
\Bigl ( \partial_{ij} \varphi(v) + (\partial_{ij}\varphi)(w) \Bigr ) \, dv dw\, ds $$
$$ \qquad\qquad {} +   \sum_{i=1}^3  \int_{t_1}^{t_2} \iint_{\R^N\times\R^N} f(s,v)\, f(s,w)\, b_i(v-w)\, 
\Bigl ( \partial_i \varphi(v) - (\partial_i \varphi)(w) \Bigr )\, dv dw\, ds. $$

 We recall that (for $i,j=1,..,3$) $a_{ij}$ and $b_i$ are defined by (\ref{annee}) and
 (\ref{coeff1}), and that thanks to estimate (\ref{cutie}), Lemma \ref{lemm2}, 
and a direct extension of the proof of Corollary~\ref{weakso}, 
each term in (\ref{starpri}) is well defined.
 \medskip
 
 Finally, we assume that $(\gamma_2+2)\,(1- \inf(\gamma_1/2, -1)) > -4$. 
 \smallskip
 
 Then, for all $k\in\R$ such that $\int_{\R^N} f_{in}(v)\, (1+ |v|^2)^k\,dv < +\infty$, the moment of order $2k$ of $f$ is locally (in time) bounded: for all $T>0$, 
 \begin{equation}\label{momct}
 \sup_{t\in [0,T]} \int_{\R^N} f(t,v)\, (1+ |v|^2)^k \, dv < +\infty . 
 \end{equation}
\smallskip

More precisely, this last quantity only depends on $T$, $N$, the parameters of $\psi$, the initial mass,
momentum, energy and (an upper bound of the) entropy, the quantity ${\mathcal{Q}}_{T,N,\gamma_1}(f)$, and
the initial moment $\int_{\R^N} f_{in}(v)\, (1+ |v|^2)^k\,dv$.
 \smallskip
 
Note that when $\gamma_1=\gamma_2$, the condition  $(\gamma_2+2)\,(1- \inf(\gamma_1/2, -1)) > -4$ is equivalent to $\gamma_1 > -2\,\sqrt{3}$, so that the Coulomb case,
 which corresponds to $\gamma_1=\gamma_2=-3$, falls within the range of application
of this proposition.
 \end{Prop}
 
 {\bf{Remarks}}: Note first that it is possible to take $\delta =0$ in this proposition, but the proof has 
then to be somewhat modified (and approximate solutions must be introduced),
 so we discarded  this case (which would be very far from any $\psi$ coming out physics anyway). 
For $\delta <0$, the behavior of moments can become completely different 
since $\psi$ can in a such a case look like a hard potential, 
we do not therefore investigate such functions $\psi$.
 \smallskip
 
As stated in the remarks of the previous section, all reasonable solutions of the Landau equation (for $\psi$ satisfying the estimates of Proposition \ref{mom1}) with initial data having a finite mass, energy and entropy, will satisfy the estimate (for all $T>0$) $\int_0^T\int_{\R^N} |\nabla \sqrt{f(t,v)}|^2\, (1+ |v|^2)^{\inf(\gamma_1/2, -1)} \, dv dt < + \infty$, thanks to Theorem \ref{princ1}, and therefore
estimate (\ref{cutie}), thanks to Lemma \ref{lemm1}. The assumptions of Proposition \ref{mom1}
are therefore quite reasonable.
 \smallskip
 
 Surprisingly, the critical parameter $\gamma = -2\,\sqrt{3}$ appears in the issue of propagation of moments. We are not sure that this critical parameter is really significant since our proof
might be far from optimal (our main purpose was to treat the Coulomb case).
\smallskip

One can observe from the proof that the dependence w.r.t. $T$ in 
estimate (\ref{momct}) is of polynomial type, since no Gronwall-type argument is used. 
\smallskip

Note that the best result (in terms of dependence of the constants w.r.t. the order of the moment) of
propagation of moments in (moderately) soft potentials is to be found in \cite{toscani_villani}.
We insist that the result of propagation of moments which is presented here is typical of
soft potentials: as for (moderately) soft potentials, no appearance of moments is expected in such a case. Propagation of (reasonable)  higher than polynomial moments can probably be obtained thanks to the summation of polynomial moments. 
\smallskip

We end up this series of remarks by the following (rather vague) statement: in some sense, the propagation of moments of any order shows that the main problems in the theory of the Landau equation 
in the Coulomb case are rather related to smoothness (or integrability) issues than to issues related to large $|v|$. 
\bigskip

The proof that we propose follows the usual lines of proofs of propagation of moments. The main
modification w.r.t. (moderately) soft potentials  consists in using the weighted $L^1_t(L^p_v)$ ($p>1$) estimate on $f$ (coming out of the 
weighted $L^2_t(H^1_v)$ norm of $\sqrt{f}$) in order to treat the singularity 
which is typical of very soft potentials.
 It is important here that the weights in the $L^1_t(L^p_v)$ estimate be polynomial,
 and not arbitrary. This was 
not the case in the previous section, where the fact that $H$-solutions are standard weak solutions only requires the information that some  $L^1_t(L^p_{v, \,loc})$ bound holds.
\smallskip

 We choose to present the proof as an induction in order to avoid an interpolation leading to a differential
 inequality that would necessitate to work on an approximate solution of the Landau equation rather than
on the (weak) solution itself. As a consequence, the proof of  Proposition \ref{mom1} is obtained by applying inductively the following lemma: 
 
 \begin{Lem}\label{emind}
 Let $N$, $T$, $\psi$ (together with $K_1$, $K_2$, $K_3$, $\delta$, $\gamma_1$, $\gamma_2$),
 and $f$ be as in Proposition \ref{mom1}.
\smallskip

 We define for all $l \in\R$,
\begin{equation}\label{ml}
 M_l(f; T) := ||f||_{L^{\infty}([0,T]; L^1_{2l} (\R^N))} = 
\sup_{t\in [0,T]} \int_{\R^N} f(t,v) \, (1+ |v|^2)^l \, dv. 
 \end{equation}
\smallskip

Let $k>1$ be such that $M_{k}(f; 0) < +\infty$, $M_{k-\delta/2}(f; T) < +\infty$.
If $\gamma_2<-2$, we assume moreover that  
 $$ M_{\frac{2\,(k-1)+(\gamma_2+2-\var)\,\inf(\gamma_1/2, -1)}{6 + \gamma_2 -\var}}(f,T) < +\infty, $$
 for some $\var\in ]0,\gamma_2 + 4[$.
\smallskip

Then, 
$$M_{k}(f; T) < +\infty. $$
More precisely, $M_{k}(f; T)$ is bounded by a constant only depending on 
 $N$, $T$, $k$, the parameters of $\psi$, the bound ${\mathcal{Q}}_{T,N,\gamma_1}(f)$, the initial moment $M_{k}(f; 0)$,  the moment $M_{k-\delta/2}(f; T)$, and, when $\gamma_2 \in ]-4, -2[$, the moment
$ M_{\frac{2\,(k-1)+(\gamma_2+2-\var)\, \inf(\gamma_1/2, -1)}{6 + \gamma_2 -\var}}(f,T)$
 and $\var\in ]0,\gamma_2 + 4[$.
 \end{Lem}
 
 {\bf{Proof of Lemma \ref{emind}}}:
 We consider $k>1$, $\chi \in D(\R)$ such that $\chi|_{[0,1]} = 1$,
 $\chi|_{[0,2]^c} = 0$, $0 \le \chi \le 1$, $ \eta \in ]0, 1[$, and we use 
  the test function $v \mapsto \phi(v) := (1+|v|^2)^{k}\, \chi(\eta\, (1 + |v|^2)^{1/2})$.
\par
Then (for $i,j = 1,..,N$),
$$ \pa_i \phi(v) = 2\,k\,v_i\,(1+|v|^2)^{k-1}\, \chi(\eta\, (1 + |v|^2)^{1/2})$$ 
$$ + \,\eta\,v_i\, (1 + |v|^2)^{k - 1/2}\, \chi'(\eta\, (1 + |v|^2)^{1/2}), $$
and
 $$ \pa_{ij} \phi(v) = 2\,k\, [(1+|v|^2)\,\delta_{ij} + 2\,(k-1)\,v_i\,v_j]\, (1 + |v|^2)^{k-2}
 \, \chi(\eta\, (1 + |v|^2)^{1/2}) $$
 $$ +\, [(4\,k-1)\,v_i\, \eta\,v_j + \eta\,\delta_{ij}\,(1+|v|^2)]\, (1 + |v|^2)^{k - 3/2}\, \chi'(\eta\, (1 + |v|^2)^{1/2}) $$
 $$ + \,(\eta\,v_i)\, (\eta\,v_j)\, (1 + |v|^2)^{k-1}\, \chi''(\eta\, (1 + |v|^2)^{1/2}), $$
 so that 
 $$ |\pa_{ij} \phi(v)| \le C(k)\,  (1 + |v|^2)^{k-1},$$
 where $C(k)$ is a constant depending only on $k$ 
(and a bound on $||\chi'||_{\infty}$, $||\chi''||_{\infty}$).
\medskip

Using definitions (\ref{annee}) and (\ref{coeff1}), we see that (for $i,j = 1,..,N$)
$$ |a_{ij}(z)| \le \psi(|z|), \qquad |b_i(z)| \le (N-1)\, |z|^{-1}\, \psi(|z|), $$
so that 
$$ | a_{ij}(v-w) \, (\pa_{ij} \phi(v) + \pa_{ij} \phi(w))  | \le C(k) \,  \psi(|v-w|)\, (1+ |v|^2 + |w|^2)^{k-1}, $$
$$ | b_i(v-w)\, (\pa_i \phi(v) - \pa_i \phi(w))|  \le C(N,k) \,  \psi(|v-w|)\, (1+ |v|^2 + |w|^2)^{k-1}, $$
where $C(k)$ (resp. $C(N,k)$) is a constant depending only on $k$ (resp. on $N,k$).
\medskip

Then, we get
 $$\int f(T,v) \,  (1+|v|^2)^{k}\, \chi(\eta\, (1 + |v|^2)^{1/2})\, dv 
 - \int f(0,v) \,  (1+|v|^2)^{k}\, \chi(\eta\, (1 + |v|^2)^{1/2})\, dv$$
$$ \le C(N,k)\, \int_0^T \int\int f(t,v)\,f(t,w)\, \psi(|v-w|)\, (1+ |v|^2 + |w|^2)^{k-1} \, dv dw dt $$
  $$ \le C(N,k) \int_0^T \int\int_{|v - w| \le 1, |v|\ge 2, |w| \ge 1}  f(t,v)\,f(t,w)\, (K_1 + K_2 \,|v-w|^{\gamma_2+2})$$
$$ \times\, (1+ |v|^2 + |w|^2)^{k-1} \, dv dw dt $$
$$ + \, C(N,k) \int_0^T \int\int_{|v - w| \le 1, |v|\le 2, |w| \le 3}  f(t,v)\,f(t,w)\, (K_1 + K_2 \,|v-w|^{\gamma_2+2})$$
$$ \times (1+ |v|^2 + |w|^2)^{k-1} \, dv dw dt $$
$$ + \, C(N,k) \int_0^T \int\int_{|v - w| \ge 1}  f(t,v)\,f(t,w)\, (K_1+ K_2)\, |v-w|^{2-\delta}$$
$$ \times\, (1+ |v|^2 + |w|^2)^{k-1} \, dv dw dt . $$
\smallskip

If  $|v- w| \le 1$ and  $|v| \ge 2$, then $ \frac12\,|v| \le |w| \le \frac32\,|v|$, so that (still when $|v- w| \le 1$ and  $|v| \ge 2$), we get the estimates
$$ 1 + |v|^2 \le 1+ |v|^2 + |w|^2 \le \frac{13}4\,(1 + |v|^2),$$
and
$$ 1 + |w|^2 \le 1+ |v|^2 + |w|^2 \le 5\,(1 + |w|^2).$$
\smallskip

As a consequence, for any $q \in \R$, 
 $$\int f(T,v) \,  (1+|v|^2)^{k}\, \chi(\eta\, (1 + |v|^2)^{1/2})\, dv
  - \int f(0,v) \,  (1+|v|^2)^{k}\, \chi(\eta\, (1 + |v|^2)^{1/2})\, dv$$
  $$ \le C(N,k)\,K_1 \, T\, M_0(f;T)\,M_{k-1}(f; T) + C(N,k)\,K_2\, \sup\left(1, \left(\frac{13}4\right)^{q/2}\right)\, \sup(1, 5^{k-1-q/2}) $$
$$ \times \,\int_0^T \int\int_{|v -w| \le 1, |v|\ge 2, |w| \ge 1}  f(t,v)\,f(t,w)\, |v-w|^{\gamma_2 +2}
\, (1 + |w|^2)^{k -1- q/2} \, (1 + |v|^2)^{q/2}\, dv dw dt $$
$$ + \, C(N,k) \,K_1 \, T\, M_0(f; T)^2$$
$$ + \,C(N,k)\,K_2 \, \int_0^T \int\int_{|v- w| \le 1, |v|\le 2, |w| \le 3}  f(t,v)\,f(t,w)\, |v-w|^{\gamma_2+2}\, dv dw dt $$
$$  +\,C(N,k)\, (K_1 + K_2)\, T \,  M_{k - \delta/2}(f; T)^2 .$$
\smallskip

If $\gamma_2 \in [-2, -\delta[$, we take $q=0$ in the estimate above, and obtain that
$$\int f(T,v) \,  (1+|v|^2)^{k}\, \chi(\eta\, (1 + |v|^2)^{1/2})\, dv \le C $$
where
$$ C:= C(T,K_1,K_2,\gamma_2, M_1(f;0), M_{k}(f; 0), M_{k - \delta/2}(f; T), N,k) $$
only depends on $N,T$, the parameters of $\psi$, the initial mass and energy, the initial moment
$M_{k}(f; 0)$, and the moment $M_{k - \delta/2}(f; T)$. 
\medskip

We now suppose that $\gamma_2 \in ]-4, -2[$. Thanks to Lemma \ref{lemm2}, 
$$ \int_0^T \int\int_{|v|\le 2, |w| \le 3}  f(t,v)\,f(t,w)\, |v-w|^{\gamma_2+2}\, dv dw dt $$
$$ \le \int_0^T  C(N,\gamma_2)\, ||f(t,\cdot)||_{L^1(\R^N)} 
\left( ||f(t,\cdot)||_{L^1_2(\R^N)} +  ||f(t,\cdot)\, 1_{|\cdot| \le 2}||_{L^{\frac{N}{N-2}}(\R^N)} \right)\, dt  $$
$$\le  C(N,\gamma_2)\, M_{0}(f; T) \, \left( M_{1}(f; T) \,T + 5^{-inf(\gamma_1/2, -1)}\, 
{\mathcal{Q}}_{T,N,\gamma_1}(f) \right) , $$ 
where $C(N,\gamma_2)$ is a constant only depending on $N$ and $\gamma_2$.
If $N=2$, one must replace $N/(N-2)$  by some $q\ge 1$ (depending on $\gamma_2$).
 \medskip
 
As a consequence, for all $r \in [1, N/(-\gamma_2-2)[$ (note that $-\gamma_2-2 < 2 \le N$) and all $q\in \R$ (and defining $r'$ by the identity $\frac1r + \frac1{r'} =1$), using Young's inequality
for convolutions, and denoting in the sequel $C(a,b,...)$ a constant depending only on $a,b, ..$,
$$ \int f(T,v) \,  (1+|v|^2)^{k}\, \chi(\eta\, (1 + |v|^2)^{1/2})\, dv $$
$$\le C(N,T,K_1,K_2,\gamma_2, \gamma_1, M_{k}(f; 0), M_{k - \delta/2}(f; T), {\mathcal{Q}}_{T,N,\gamma_1}(f), q,k)\,$$
$$ \times
\bigg( 1 + \int_0^T  ||(1+ |\cdot|^2)^{q/2}\, f||_{L^{r'}(\R^N)} $$
$$ \times \,  ||(1+ |\cdot|^2)^{k-1-q/2}\, f||_{L^{1}(\R^N)}
\, || z \mapsto |z|^{\gamma_2+2} 1_{|z| \le 1} ||_{L^r(\R^N)} \, dt \bigg) $$
$$\le C(N,T,K_1,K_2,\gamma_2, \gamma_1, M_{k}(f; 0), M_{k - \delta/2}(f; T), {\mathcal{Q}}_{T,N,\gamma_1}(f), q,k,r)\,$$
$$ \times
\bigg( 1 + \int_0^T  ||(1+ |\cdot|^2)^{q/2}\, f||_{L^{r'}(\R^N)} \,  ||(1+ |\cdot|^2)^{k-1-q/2}\, f||_{L^{1}(\R^N)}
 \, dt \bigg) .$$
Taking (for some $\var\in ]0, \gamma_2 + 4[$) $r' = \frac{N}{N + \gamma_2 +2-\var}$, $\beta = \frac{-\gamma_2 -2+\var}2 \in ]0,1[$, 
and $q = \frac{\beta\, \inf(\gamma_1,-2) + 2\,(k-1)\,(1-\beta)}{2-\beta}$,
we see that $r \in [1, N/(-\gamma_2-2)[$, and
$$ \frac1{r'} = \beta\, \frac{N-2}{N} + 1-\beta, \qquad  q = \beta\, \inf(\gamma_1,-2) + (1-\beta)\,(2\,(k-1)-q),$$ 
so that thanks to the interpolation inequality of Proposition \ref{interp2} in Appendix 1,
$$ ||(1+ |\cdot|^2)^{q/2}\, f||_{L^{r'}(\R^N)} \le ||(1+ |\cdot|^2)^{\inf(\gamma_1/2,-1)} f||_{L^{\frac{N}{N-2}} (\R^N)}^{\beta} \,
||(1+ |\cdot|^2)^{k-1-q/2} f||_{L^{1} (\R^N)}^{1-\beta}. $$  
 We end up with
$$ \int f(T,v) \,  (1+|v|^2)^{k}\, \chi(\eta\, (1 + |v|^2)^{1/2})\, dv$$
$$ \le C(N,T,K_1,K_2,\gamma_2, M_{k}(f; 0), M_{k-\delta/2} (f; T), {\mathcal{Q}}_{T,N,\gamma_1}(f), \var, k,\gamma_1)$$
$$ \times\, \bigg( 1 + \int_0^T  ||(1+ |\cdot|^2)^{\inf(\gamma_1/2, -1)} f||_{L^{\frac{N}{N-2}}}^{\beta} \,
  ||(1+ |\cdot|^2)^{k-1-q/2} f||_{L^{1}}^{2-\beta}\, dt \bigg)$$
  $$ \le C \bigg(N,T,K_1,K_2,\gamma_2, M_{k}(f; 0), M_{k-\delta/2}(f; T),$$
$$ M_{\frac{2\,(k-1)+(\gamma_2+2-\var)\,\inf(\gamma_1/2,-1)}{6 + \gamma_2 -\var}}(f,T), {\mathcal{Q}}_{T,N,\gamma_1}(f),\var,k,\gamma_1 \bigg) , $$
since $k-1 - q/2 = \frac{2\,(k-1)+(\gamma_2+2-\var)\,\inf(\gamma_1/2,-1)}{6 + \gamma_2 -\var}$.
\smallskip

We conclude the proof of Lemma \ref{emind} by letting $\eta$ tend to $0$ and by using Fatou's lemma.
$\square$
\bigskip

We now turn to the
\smallskip

{\bf{Proof of Proposition \ref{mom1}}}: Since we assume that $(\gamma_2+2)\,(1 - \inf(\gamma_1/2,-1)) > -4$, we can find $\var \in ]0, \gamma_2+4[$
such that $(\gamma_2+2 -\var)\,(1 - \inf(\gamma_1/2,-1)) > - 4$. Then we define by induction (for $n \in \N$) the sequence
$$ k_0 = k, \qquad k_{n+1} = \sup \bigg( k_n - \delta/2,  \frac{2}{6+\gamma_2 -\var} \,k_n+ \frac{ \inf(\gamma_1/2,-1)\,(\gamma_2+2-\var) -2}{6 + \gamma_2 - \var}\bigg),$$
when $\gamma_2 \in ]-4, -2[$, and
$$ k_0 = k, \qquad k_{n+1} =  k_n - \delta/2, $$
when $\gamma_2\in [-2, -\delta]$.
\smallskip

We observe that the fixed point of $k \mapsto \frac{2}{6+\gamma_2 -\var} \,k+ \frac{ \inf(\gamma_1/2,-1)\,(\gamma_2+2-\var) -2}{6 + \gamma_2 - \var}$ is
$\bar{k} =  \frac{ \inf(\gamma_1/2,-1)\,(\gamma_2+2-\var) -2}{4 + \gamma_2 - \var}$, 
so that $ \bar{k} < 1$ is equivalent to $(\gamma_2+2 -\var)\,(1 - \inf(\gamma_1/2,-1)) > -4$, which holds
thanks to our assumption and our choice of $\var$.
Then, the sequence $(k_n)_{n\in\N}$ decreases, and we can pick up $n_0\in\N$ such that $k_{n_0} > 1$, and $k_{n_0 +1} \le 1$. 
 By (finite) induction, starting from $k_{n_0}$ and
using Lemma~\ref{emind} repeatedly, we get
 that $M_{k_n}(f,T) < +\infty$ for all $n = n_0,n_0-1,..,0$, so that
Proposition \ref{mom1} holds.
$\square$
\medskip

{\bf{Proof of Proposition \ref{momco}}}: We first observe that $\psi(z) = |z|^{-3}$ satisfies
the assumption on $\psi$ in Proposition \ref{mom1}, with $\gamma_1=\gamma_2 = -3$, so that
$(\gamma_2+2)\,(1 - \inf(\gamma_1/2,-1)) = - 5/2 >-4$. Then, the H-solutions of the Landau equation in
the Coulomb case satisfy the conservation of mass, momentum and energy, and have an upper-bounded
entropy. According to Corollary \ref{weakso}, they also lie in  $L^{\infty}([0,T]; L^3_{-3}(\R^3))$ for all $T>0$, and (cf. \cite{vill:new:97}) satisfy the (strong w.r.t. $t$) weak form of the Landau equation. Proposition~\ref{momco} is then a direct consequence of Proposition~\ref{mom1}. 
$\square$.

\section{Conclusions and Perspectives}

 In this paper, we have shown that the existing solutions $f$ of the Landau equation in the Coulomb case (that is, the $H$-solutions) are such that $\sqrt{f}$ lies in a weighted $L^2(H^1)$
 space, so that $f$ lies in a weighted $L^1(L^3)$ space. This is sufficient to show that those 
solutions can be defined in the standard weak framework, and that all polynomial $L^1$
moments are propagated. 
 \smallskip
 
 According to \cite{kcw}, such an estimate enables the propagation of $L^p$ norms (for $p>1$)
when (moderately) soft potentials are concerned (we also refer to Appendix~2 of this work).
However, we weren't able to extend this result in the Coulomb case (or for any $\gamma < -2$),
so that the smoothness issues for the (spatially homogeneous) Landau equation in the Coulomb case
are still completely open (together with the related issue of uniqueness). 
\smallskip

We believe that the results presented here for the Landau equation in the Coulomb (or more generally  soft potentials) case are related to the results of \cite{gs1}, \cite{gs2}, \cite{gs3} for the Boltzmann equation without angular cutoff. Establishing an explicit link with 
those results seems to us a very interesting line of research. 
\bigskip
 
 \begin{center}
{\bf{Appendix 1: Interpolations}}
\end{center} 
\bigskip

We recall two classical interpolation inequalities, which can be proven by an application of H\"older's inequality.
\medskip




The first one deals with weighted $L^p$ spaces.
 
\begin{Prop} \label{interp2}
 Assume that for some $q_1, q_2 \in [1,\infty]$, $a_1, a_2 \in \R$,  $x \mapsto (1+|x|^2)^{a_1} \,f(x)\in L^{q_1}(\R^N)$,
 $x \mapsto (1+|x|^2)^{a_2} \,f(x)\in L^{q_2}(\R^N)$.
\medskip

 Then, for and any $\beta \in [0,1]$, 
$x \mapsto (1+|x|^2)^a\, f\in L^{q}(\R^N)$, where
$$ \frac1q = \frac{\beta}{q_1} + \frac{1-\beta}{q_2}, \qquad 
 a = \beta \, a_1 + (1-\beta)\,a_2. $$
Moreover,
\begin{equation}\label{interpf2}
||x \mapsto (1+|x|^2)^a\,  f(x)||_{L^q(\R^N)} \le  ||x \mapsto (1+|x|^2)^{a_1} \,f(x)||_{L^{q_1}(\R^N)}^{\beta}
\end{equation}
$$ \times\, ||x \mapsto (1+|x|^2)^{a_2}\, f(x)||_{L^{q_2}(\R^N)}^{1-\beta}, $$
or, in abridged form,
\begin{equation}\label{interpmaxbibi}
||f||_{L^q_a(\R^N) }  \le 
||f||_{L^{q_1}_{a_1}(\R^N)}^{\beta}\,\,||f||_{L^{q_2}_{a_2}(\R^N) }^{1-\beta} .
\end{equation}
\end{Prop}
\bigskip

Then, the second one
 deals with weighted $L^p(L^q)$ spaces, the weight 
 concerning the second variable. 

\begin{Prop} \label{interp3}
 Assume that for some $p_1,q_1, p_2,q_2 \in [1,\infty]$, $a_1, a_2 \in \R$,  $(t,x) \mapsto (1+|x|^2)^{a_1} \,f(t,x)\in L^{p_1}([0,T];\, L^{q_1}(\R^N))$ on one hand,
and $(t,x) \mapsto (1+|x|^2)^{a_2} \, \,f(t,x)\in L^{p_2}([0,T];\, L^{q_2}(\R^N))$
 on the other hand.
\smallskip

 Then, for any $\beta \in [0,1]$, 
$$(t,x) \mapsto (1+|x|^2)^{a} \, f(t,x)\in L^{p}([0,T]; \, L^{q}(\R^N)),$$
 where
$$ \frac1p = \frac{\beta}{p_1} + \frac{1-\beta}{p_2}, \qquad 
 \frac1q = \frac{\beta}{q_1} + \frac{1-\beta}{q_2}, \qquad a = \beta \, a_1 + (1-\beta)\,a_2. $$
Moreover,
\begin{equation}\label{interpf1}
||(t,x) \mapsto (1+|x|^2)^{a} \, f(t,x)||_{L^p([0,T];\, L^q(\R^N)) } 
\end{equation}
$$\le  ||(t,x) \mapsto (1+|x|^2)^{a_1} \,f(t,x)||_{L^{p_1}([0,T];\, L^{q_1}(\R^N)) }^{\beta}$$
$$ \times\,  ||(t,x) \mapsto (1+|x|^2)^{a_2} \,f(t,x)||_{L^{p_2}([0,T];\, L^{q_2}(\R^N)) }^{1-\beta}, $$
or, in abridged notation,
\begin{equation}\label{interpmax}
||f||_{L^p([0,T]; \,L^q_a(\R^N)) }  \le 
||f||_{L^{p_1}([0,T];\, L^{q_1}_{a_1}(\R^N)) }^{\beta}\,\,||f||_{L^{p_2}([0,T];\, L^{q_2}_{a_2}(\R^N)) }^{1-\beta} .
\end{equation}
\end{Prop}
\bigskip
 
 \begin{center}
{\bf{Appendix 2: Further regularity estimates in the case of (moderately) soft potentials}}
\end{center} 
\bigskip


 This second appendix is devoted to the application of the main estimate of the paper (that is, Theorem 
 \ref{princ1}), to the case of moderately soft potentials. Since a well-established theory exists in this
 case (cf. in particular \cite{kcw} and \cite{alex_new}), we merely present alternative proofs of already
existing results, or proofs of variants of those results. 
\medskip 

 We consider functions $\psi$ in the Landau kernel which are bounded above around point $0$ by power laws
of (moderately) soft potentials type. We also choose to impose a bound from below involving an arbitrary power 
at infinity, and a bound above at infinity which behaves like a  (moderately) soft potential, but with an exponent which needn't be the same as the exponent appearing around point $0$.
\smallskip

In order to keep the proofs of our results as simple as possible,
we consider smooth solutions of the Landau equation rapidly decaying at infinity (w.r.t. $v$)
such that $\ln f$ is slowly growing at infinity (w.r.t. $v$),
in the presentation of the results. In other words, we only describe the {\it{a priori}} estimates. 
\par
Moreover, for the sake of readability of the results, we present estimates which depend on
an arbitrary number of moments (in $L^1$) on the initial data. Those results still (partially)
 hold when only a finite number of moments initially exist.
\smallskip

The first part of the estimates of \cite{kcw},
 that is the bound in a weighted $L^1([0,T]; L^3(\R^3))$
 (or, in dimension $N\ge 2$, in a weighted $L^1([0,T]; L^{N/(N-2)}(\R^N))$),
 and by interpolation with $L^{\infty}([0,T]; L^1(\R^3))$,
 in a weighted $L^{1+\delta}([0,T]; L^{3 - \var}(\R^3))$
for some $\delta, \var>0$,
is a direct consequence of the main estimate
of this paper and Sobolev inequalities (cf. Theorem \ref{princ1}, and Lemma \ref{lemm1}).
\smallskip

The second part of the estimates of \cite{kcw} are related to the propagation of $L^p$ norms. 
The bounds that we present in this appendix also concern the propagation of $L^p$ norms. The main 
differences with the results of \cite{kcw} are the following: 
\begin{itemize}
\item
We use the spaces $\cup_{p<q} L^q$ whereas the spaces $L^p$ are used in \cite{kcw};
\item
We assume here that all moments in $L^1$ of the solution are (initially) finite, whereas only a finite number of them are assumed
to be  finite in \cite{kcw};
\item
 Our estimates are polynomially growing w.r.t. time, whereas they grow more rapidly in \cite{kcw};
\item
 More importantly from the structural point of view, in the proof of the bounds, the coercivity estimate $ \sum_{i,j} (a_{ij} * f) \,\xi_i\,\xi_j \ge C\,(1 +|v|^2)^{\gamma/2}\,|\xi|^2$ proven in \cite{kcw}, \cite{alex_new}, is replaced by our main theorem (Theorem \ref{princ1});
\item 
We do not treat here the limiting 
 case $\gamma = -2$, whereas this case is treated in \cite{kcw};
\item
Finally, we authorize here dimensions different from $3$ and functions $\psi$ which are
 not necessarily power laws.
 \end{itemize}
\medskip

We start with the most technical case, which appears  when $\psi(z) \le K_2\, |z|^{\gamma_2+2}$ with $\gamma_2\in ]-2, -1[$ (for $|z|$ close to $0$), 
and ones looks for the propagation of $L^q$ norms with $q<1 + \frac1{N + 1 + \gamma_2}$. 
This corresponds to a soft potential in some intermediate regime (still moderately soft, but already far
from Maxwell molecules),
 and a parameter $q$ close to $1$.
\medskip

 Our proposition writes:

\begin{Prop}\label{proreg1}
 Let $\psi$ be a function satisfying
$$ \forall z\ge 0, \qquad K_3\, \inf(1, |z|^{\gamma_1 + 2}) \le \psi(z) \le  K_1\, |z|^{2-\delta} + K_2\, |z|^{\gamma_2+2}, $$
for some $K_1,K_2, K_3>0$ and $\delta \in ]0,1]$, $\gamma_1 \le 0$, $ \gamma_2\in ]-2, -1[$. 
\par
 Let $T>0$, and $f:=f(t,v)\ge 0$ be a strong solution of the Landau equation (\ref{Landau}) on $[0,T]\times\R^N$
(more precisely, we assume that $f$ is smooth, rapidly decaying when $|v| \mapsto +\infty$.
 We also assume that $\ln f$ is slowly growing when $|v| \mapsto +\infty$).
 \smallskip
 
We consider $q\in ]1, 1 + \frac1{N + 1 + \gamma_2}]$.
Then, for all $m\in\R$, $\var_0 \in ]0, q-1]$, there exists $C$ depending only on
$N,K_1,K_2,K_3,\gamma_1,\gamma_2,\delta$, $\var_0,m,q$, on the initial mass, momentum, energy and (upper part of) the entropy, on the initial moments (in $L^1$) 
$(||f(0, \cdot)||_{L^1_{\mu}(\R^N)})_{\mu >1}$, 
and initial $L^{q -\var}$ norm
 $(||f(0, \cdot)||_{L^{q-\var}(\R^N)})_{\var \in ]0, q-1]}$, such that
$$ \sup_{t\in [0,T]} ||f(t, \cdot)||_{L^{q - \var_0}_m(\R^N)} \le C. $$
 \end{Prop}
 
This proposition is obtained by using inductively the following lemma:

\begin{Lem}\label{ebody}
Let $\psi$ and $f$ satisfy the assumptions of Proposition \ref{proreg1}.
\par
We consider $A \in ]1, \frac{N}{N+1+\gamma_2}]$.
\smallskip

Then for all $\var_0 \in ]0, A-1]$ and $m_0\in\R$,
$$ \sup_{t\in [0,T]} ||f(t, \cdot)||_{L^{A\,(1 + (\gamma_2+2)/N) -\var_0}_{m_0}(\R^N)} $$
is bounded by a constant depending on the same parameters as in Proposition \ref{proreg1}
(except that $q$ is replaced by $A$, and the initial $L^{q -\var}$ norm 
by the initial $L^{A -\var}$ norm),
and on the quantities 
$$ \sup_{t\in [0,T]} ||f(t, \cdot)||_{L^{A-\var}_m(\R^N)} ,$$
for all $\var \in ]0, A-1]$, and $m\in\R$.
 \end{Lem}

{\bf{Proof of Lemma \ref{ebody}}}:
We observe that multiplying Landau's equation by $f^k$ for some $k>0$ and integrating on $[0,T] \times \R^N$,
we end up with 
$$ \int \frac{f^{k+1}}{k+1}\, dv\, (T) + k \int_0^T \int  \sum_{i=1}^N\sum_{j=1}^N f^{k-1}\, (a_{ij} * f) \, \pa_i f \, \pa_j f \, dv dt $$
$$ = \int \frac{f^{k+1}}{k+1}\, dv \,(0) + k  \,\int_0^T \int  \sum_{i=1}^N (b_i * f) \, \pa_i f \, f^k\, dv $$
$$ \le \int \frac{f^{k+1}}{k+1}\, dv \,(0) + k \,{\mathcal{F}}_{T,N,\gamma_1}(f)^{1/2} \, \bigg(\int_0^T \int  \left|\sum_{i=1}^N b_i * f\right|^2
\, f^{2k+1}\, (1+|v|^2)^{\sup(-\gamma_1/2, 1)}\, dv dt \bigg)^{1/2} ,$$
with, thanks to Theorem \ref{princ1},
$$ {\mathcal{F}}_{T,N,\gamma_1}(f) := \int_0^T \int_{\R^N} \frac{|\nabla_v f|^2}{f}  \,
 (1 + |v|^2)^{\inf(\gamma_1/2, -1)}\, dv  dt \le C, $$
 where $C$ depends only on $N$, $K_3$, $\gamma_1$, and the initial mass, momentum, energy and 
 (an upper bound of the) initial entropy. In the rest of the proof of this lemma, we denote by $C(a,b,..)$
 a constant only depending on $a,b,..$.
 \smallskip
 
Then, for all $r\in [1,\infty]$,
$$ \int_0^T \int  \left|\sum_{i=1}^N b_i * f \right|^2\, f^{2k+1}\, (1+|v|^2)^{\sup(-\gamma_1/2, 1)}\, dv dt $$
$$\le
\, C(N) \, \int_0^T \int  f^{2k+1}(t,v) \bigg| \int_{|v-w|\ge 1} (K_1\,|v-w| + K_2)\, f(t,w)\,dw \bigg|^2\, (1+|v|^2)^{\sup(-\gamma_1/2, 1)}\, dv dt $$
$$ + \, C(N) \, \int_0^T \int  f^{2k+1}(t,v) \bigg| \int_{|v-w|\le 1} (K_1 + K_2\, |v-w|^{\gamma_2+1})\, f(t,w)\,dw \bigg|^2\, (1+|v|^2)^{\sup(-\gamma_1/2, 1)}\, dv dt $$
 $$ \le C(N,K_1,K_2, ||f(0, \cdot)||_{L^1_2(\R^N)})\, \int_0^T \int  f^{2k+1}(t,v) \, (1+|v|^2)^{\sup(1-\gamma_1/2, 2)}\, dv dt $$
$$+ \,C(N,K_2, ||f(0, \cdot)||_{L^1(\R^N)})\,
 \int_0^T \int  f^{2k+1}(t,v)  \, | \,|\cdot|^{\gamma_2+1}\,1_{|\cdot|\le 1} *f(t,\cdot) |^2 (v) \, (1+|v|^2)^{\sup(-\gamma_1/2, 1)}\, dv dt$$
 $$ \le C(N,K_1,K_2, ||f(0, \cdot)||_{L^1_2(\R^N)})
 \, || f ||_{L^{2k+1}([0,T]; L^{2k+1}_{ \frac{\sup(2-\gamma_1, 4)}{2k +1} } (\R^N))}^{2k+1} $$
$$ + \,C(N,K_2, ||f(0, \cdot)||_{L^1(\R^N)} )
 \int_0^T  || (1+|\cdot|^2)^{\sup(-\gamma_1/2,1)}\,f^{2k+1}(t,\cdot)||_{L^{r'}(\R^N)}  \, 
|| \,\,|\cdot|^{\gamma_2+1}\,1_{|\cdot|\le 1} * f(t,\cdot) ||_{L^{2r}(\R^N)}^2 \, dt, $$
for all $r\in [1, +\infty]$ and $r'$ defined by $\frac1r + \frac1{r'}=1$.

Note that taking $\beta = \frac{N}{N + 2 (A -\var_0)}$ (with $\var_0 \in ]0, A-1]$), one has $\beta\in ]0,1[$ and (with a slight abuse of notation)
$$ \frac{N}{N + 2 (A -\var_0)} = \frac{1-\beta}{\infty} + \frac{\beta}{1} = \frac{1-\beta}{A -\var_0} + \beta\, \frac{N-2}N, $$
when $N\ge 3$, so that, thanks to the interpolation described in Proposition \ref{interp3}, 
$$ ||f||_{L^{\frac{N + 2 (A -\var_0)}N} ([0,T];\, L^{\frac{N + 2 (A -\var_0)}N}_{m_0}(\R^N))} $$ 
is bounded by a constant which depends on 
$||f||_{L^{\infty}([0,T];\, L^{A -\var_0}_m(\R^N))}$ with some $m\in \R$
 (this quantity appears in the assumption of the
lemma) and $||f||_{L^{1}([0,T];\, L^{\frac{N}{N-2}}_{\inf(\gamma_1, -2)}(\R^N))}$, which is itself
bounded by a constant depending on $ {\mathcal{F}}_{T,N,\gamma_1}(f)$, thanks to Lemma~\ref{lemm1}.
Then, $ || f ||_{L^{2k+1}([0,T];\, L^{2k+1}_{ \frac{\sup(2-\gamma_1,4)}{2k+1}} (\R^N))}$
is bounded (by a constant depending on the parameters
 appearing in the assumptions of the lemma) for all $k < A/N$. The case
$N=2$ leads to the same conclusion if one replaces $\frac{N}{N-2}$ by any $q \in [1, +\infty[$.
\medskip

We now take $r$ in such a way that 
$$\frac{1}{2r} = \frac1{A-\var_0} - \frac{\gamma_2+1}N - 1 
               = \frac1{A-\var_0/2} - \frac{\gamma_2+1}N\,\left[ 1 - \frac{N}{\gamma_2+1}\,\frac{\var_0/2}{(A-\var_0)\,(A - \var_0/2)} \right] - 1 , $$
 (with $\var_0 \in ]0, A-1[$ small enough). Note that $-(\gamma_2+1) <N/2$ since $N\ge 2$. Moreover $\frac1{A-\var_0} - \frac{\gamma_2+1}N - 1 >0$ since $A\le \frac{N}{N +\gamma_2+1}$. 
Then, $r>1$, and thanks to Young's 
 inequality for convolutions, since $ |\cdot|^{\gamma_2+1}\,1_{|\cdot|\le 1} \in L^{ \frac{N}{-\gamma_2-1} - \delta}$
 for all $\delta >0$ small enough,
$$ \sup_{t\in [0,T]} \, || \,|\cdot|^{\gamma_2+1}\,1_{|\cdot|\le 1} * f(t,\cdot) ||_{L^{2r}(\R^N)}$$
is bounded by a constant depending only on $||f||_{L^{\infty}([0,T]; L^{A -\var_0/2}(\R^N))}$. 
\medskip

It remains to bound 
$$ \int_0^T  ||v \mapsto (1+|v|^2)^{\sup(-\gamma_1/2, 1)}\,f^{2k+1}(t,v)||_{L^{r'}(\R^N)}  \, dt . $$
 Note that taking $\beta = \frac1{2\,(A-\var_0)\,(1+ \frac{2+\gamma_2}N) - 1}$, one has $\beta \in ]0,1[$ (remember that 
 $A -\var_0 >1$), and (with a slight abuse of notation)
 $$  \frac1{2\,(A-\var_0)\,(1+ \frac{2+\gamma_2}N) - 1} = \frac{1-\beta}{\infty} + \frac{\beta}{1},$$
$$ \frac1{2\,(A-\var_0)\,(1+ \frac{2+\gamma_2}N) - 1} \, \frac1{r'} = \frac{1-\beta}{A -\var_0} + \beta\, \frac{N-2}N, $$
so that 
$$ ||f||_{L^{2\,(A-\var_0)\,(1+ \frac{2+\gamma_2}N) - 1}([0,T]; \, L^{[2\,(A-\var_0)\,(1+ \frac{2+\gamma_2}N) - 1]\,r'}_{\frac{\sup(-\gamma_1, 2)}{2k+1}} (\R^N)) }$$
is bounded by a constant depending only on 
$ ||f||_{L^{\infty}([0,T];\, L^{A -\var_0}_m(\R^N))}$  for all $m\in\R$ 
(which is controlled thanks to the assumptions of the lemma) and  $||f||_{L^{1}([0,T]; \, L^{\frac{N}{N-2}}_{\inf(\gamma_1, -2)}(\R^N))}$, which is itself
bounded by a constant depending on $ {\mathcal{F}}_{T,N,\gamma_1}(f)$, thanks to Lemma~\ref{lemm1}.
Then , for  $k< A\,(1+ \frac{2+\gamma_2}N) - 1$, we see that
$2k+1\le  2\,(A-\var_0)\,(1+ \frac{2+\gamma_2}N) - 1$ for some $\var_0 \in ]0, A-1[$ small enough),
 so that for such a value of the parameter $k$,
 $ \int_0^T  ||v \mapsto (1+|v|^2)^{\sup(-\gamma_1/2, 1)}\,f^{2k+1}(t,v)||_{L^{r'}(\R^N)}  \, dt$
 can be bounded by the quantities appearing in the statement of the lemma.
 \medskip
 Finally, we observe that thanks to Proposition \ref{mom1} (and Lemma \ref{lemm1}), $f$ has bounds in $L^{\infty}([0,T]; L^1_M(\R^N))$ for any $M\in \R$, depending only on  quantities appearing in the statement of the lemma. Thanks to a last interpolation, we conclude the proof of Lemma \ref{ebody}.
 $\square$
 \medskip
 
 We now turn to the
\medskip

{\bf{Proof of Proposition \ref{proreg1}}}: We fix $M\in \N$ as the smallest integer such that $q\,(1+ \frac{\gamma_2+2}N)^{-M} \le 1$. Since $q \le 1 + \frac1{N + 1 + \gamma_2}$, we see that $q\,(1+ \frac{\gamma_2+2}N)^{-1} \le \frac{N}{N+1+\gamma_2}$. 
Therefore we can apply Lemma \ref{ebody} by induction for $A = 1$ [in this first step, one uses a slight variant of the lemma, which holds in theory only for $A>1$],
 $A=  q\,(1+ \frac{\gamma_2+2}N)^{-M+1}$, .. , $A=q\,(1+ \frac{\gamma_2+2}N)^{-1}$. We end up with the conclusion of Proposition 
\ref{proreg1}.
 $\square$
 \bigskip
 
 We now consider a less technical case, which appears  either when $\psi(z) \le K_2\, |z|^{\gamma_2+2}$ with $ \gamma_2\in ]-2, -1[$ (for $|z|$ close to $0$), and when one looks for the propagation
of $L^q$ norms with $q\in [1 + \frac1{N + 1 + \gamma_2}, \frac{N}{N-1}[ $, or when 
 $q\in ]1, \frac{N}{N-1}[ $ and
the soft potential is close to Maxwellian molecules, that is $ \gamma_2\in [-1,-0[$.
The specificity of this case is that in the proof of the corresponding proposition (Proposition \ref{proreg2}),
 $b_i * f$ lies in $L^{\infty}_{loc} ([0,T] \times  \R^N)$, which leads to great simplifications in the proof (with respect to the proof
of Proposition~\ref{proreg1}).
\medskip

 Our proposition writes
 
 \begin{Prop}\label{proreg2}
 Let $\psi$ satisfying
$$ \forall z\ge 0, \qquad K_3\, \inf(1, |z|^{\gamma_1 + 2}) \le \psi(z) \le  K_1\, |z|^{2-\delta} + K_2\, |z|^{\gamma_2+2}, $$
for some $K_1,K_2, K_3>0$ and $\delta \in ]0,2[$, $\gamma_1 \le 0$, $ \gamma_2\in ]-2, -\delta[$. 
\par
 Let $T>0$, and $f:=f(t,v)\ge 0$ be a strong solution of the Landau equation (\ref{Landau}) on $[0,T]\times\R^N$
(more precisely, we assume that $f$ is smooth, rapidly decaying when $|v| \mapsto +\infty$.
 We also assume that $\ln f$ is slowly growing when $|v| \mapsto +\infty$).
 \smallskip

We consider $q\in [1 + \frac1{N + 1 + \gamma_2}, \frac{N}{N-1}[ $ if $ \gamma_2\in ]-2,-1[$, 
and $q\in ]1, \frac{N}{N-1}[$ if $\gamma_2\in [-1,-0[$.
\smallskip

Then, for all $m\in\R$, $\var_0 \in ]0, q-1]$, there exists $C$ depending only on
$N,K_1,K_2,K_3,\gamma_1,\gamma_2,\delta$, $\var_0,m, q$, on the initial mass, momentum, energy and (upper part of) the entropy, on the initial moments (in $L^1$) 
$(||f(0, \cdot)||_{L^1_{\mu}(\R^N)})_{\mu >1}$, 
and initial $L^{q -\var}$ norm
 $(||f(0, \cdot)||_{L^{q-\var}(\R^N)})_{\var \in ]0, q-1]}$, such that
$$ \sup_{t\in [0,T]} ||f(t, \cdot)||_{L^{q - \var_0}_m(\R^N)} \le C. $$

 

 \end{Prop}

 Once again, the proposition is obtained thanks to the inductive use of a technical lemma:

\begin{Lem}\label{ebody2}
Let $\psi$ and $f$ satisfy the assumptions of Proposition \ref{proreg2}.
\par
We consider $A \in ]\frac{N}{N + 1 + \gamma_2}, \frac{N}{N-1}[$ if $ \gamma_2\in ]-2,-1[$ and 
$A\in [1, \frac{N}{N-1}[$ if $ \gamma_2\in [-1,-0[$.
\smallskip

Then for all $\var_0 \in ]0, A-1]$ and $m_0\in\R$,
$$ \sup_{t\in [0,T]} ||f(t, \cdot)||_{L^{1 + (A -\var_0)/N }_{m_0}(\R^N)} $$
is bounded by a constant depending on the same parameters as in Proposition \ref{proreg1} (except that $q$ is replaced by $A$, and the initial $L^{q -\var}$ norm 
by the initial $L^{A -\var}$ norm), and on the quantities 
$$ \sup_{t\in [0,T]} ||f(t, \cdot)||_{L^{A-\var}_m(\R^N)} ,$$
for all $\var \in ]0, A-1]$, and $m\in\R$.
 \end{Lem}

{\bf{Proof of Lemma \ref{ebody2}}}: As in Lemma \ref{ebody}, we have (for any $k>0$)
\begin{equation}\label{uut}
 \int \frac{f^{k+1}}{k+1}\, dv\, (T) + k \int_0^T \int f^{k-1}\, \sum_{i=1}^N\sum_{j=1}^N (a_{ij} * f) \, \pa_i f \, \pa_j f \, dv dt 
 \end{equation}
$$ \le \int \frac{f^{k+1}}{k+1}\, dv \,(0) $$
$$+\, C(k,T,N,\gamma_1)\, \bigg(\int_0^T \int  \left|\sum_{i=1}^N b_i * f \right|^2
\, f^{2k+1}\, (1+|v|^2)^{\sup(-\gamma_1/2, 1)}\, dv dt \bigg)^{1/2} ,$$
where $C(a,b,..)$ is (here and in the rest of the lemma) a constant depending only on $a,b,..$.
\smallskip

We compute (for $i=1,..,3$)
$$ | b_i * f(t,\cdot)|(v) \le C(N,K_1,K_2)\, \bigg( |\,|\cdot|^{1-\delta} * f(t,\cdot)|(v)  
+ |\,|\cdot|^{1+\gamma_2} * f(t,\cdot)|(v) \bigg) $$
$$ \le C(N,K_1,K_2, ||f(t,\cdot)||_{L^1_1(\R^N)}) \,
\left(1 + |v| + |\,|\cdot|^{1+\gamma_2}\,1_{|\cdot|\le 1} * f(t,\cdot)|(v) \right). $$
Then, we use Young's inequality for convolutions. If $\gamma_2 \in [-1,0[$, we directly obtain
that $\sup_{t\in [0,T]} || b_i * f(t,\cdot)||_{L^{\infty}_{-1}(\R^N)}$ is bounded by a constant 
only depending on the parameters appearing in the statement of the lemma. If $\gamma_2 \in ]-2,-1[$,
we use the bound $ \sup_{t\in [0,T]} ||f(t, \cdot)||_{L^{A-\var}(\R^N)},$
where $A > \frac{N}{N + 1 + \gamma_2}$, and $\var>0$ is sufficiently small. We also get that
$\sup_{t\in [0,T]} || b_i * f(t,\cdot)||_{L^{\infty}_{-1}(\R^N)}$ is bounded by a constant 
only depending on the parameters appearing in the statement of the lemma.
\medskip

Then, it remains to bound
 $\int_0^T \int  f^{2k+1}\, (1+|v|^2)^{\sup(1-\gamma_1/2,2)}\, dv dt $,
with $k = \frac{A-\var_0}N$.
\par
We consider  $\beta = \frac{1}{2k+1} = \frac{N}{N+2\,(A-\var_0)}$, so that
$\beta\in ]0,1[$ and (with a slight abuse of notation)
$$ \frac{1}{2k+1} = \frac{1-\beta}{\infty} + \frac{\beta}{1} 
= \frac{1-\beta}{A -\var_0} + \beta\, \frac{N-2}N. $$
Using the bound assumed in the lemma 
(more precisely, $||f||_{L^{\infty}([0,T];\, L^{A -\var_0}_m(\R^N))}$ for all $m\in\R$),
 and remembering
that $||f||_{L^{1}([0,T]; \, L^{\frac{N}{N-2}}_{\inf(\gamma_1,-2)}(\R^N))}$ 
is bounded thanks to Theorem \ref{princ1} and Lemma \ref{lemm1} (in dimension $2$, the quantity $N/(N-2)$ is replaced by all 
$q \ge 1$), we see that the interpolation described in Proposition 
\ref{interp3} implies that
$ || f ||_{L^{2k+1}([0,T]; L^{2k+1}_{\frac{\sup(2-\gamma_1,4)}{2k+1}} (\R^N))}$ is bounded by a constant 
only depending on the parameters in the statement of the lemma.
\par
A last interpolation, like at the end of Lemma \ref{ebody}, leads to the conclusion
 of the proof of Lemma \ref{ebody2}.
 $\square$
 \medskip
 
 We now turn to the
\medskip

{\bf{Proof of Proposition \ref{proreg2}}}: We build the sequence $q_0 := q$, and 
$q_{n+1} := (q_n - 1)\,N$, for all $n\in\N$.
We fix $M\in \N$ as the smallest integer such that $q_M\le 1$ if $\gamma_2\in  [-1,0[$,
 and $q_M \le 1 + \frac{1}{N+1+\gamma_2}$ if $\gamma_2\in  ]-2,-1[$. 
Since $q < \frac{N}{N - 1}$, we see that such an integer exists.
\medskip

If $\gamma_2\in  [-1,0[$, we can apply Lemma \ref{ebody2} by induction for $A = 1$ 
(for this step, one uses a slight variant of the lemma),  $A= q_{M-1}$, .. , $A=q_1$, and get
the conclusion of Proposition \ref{proreg2}.
\medskip

If $\gamma_2\in  ]-2,-1[$, then $\frac{N}{N+1+\gamma_2} < q_M  \le 1 + \frac{1}{N+1+\gamma_2}$.  Using the result of Proposition \ref{proreg1}, we can apply Lemma \ref{ebody2} by induction for 
$A = q_M$, $A= q_{M-1}$, .. , $A= q_1$, and also get the conclusion of Proposition \ref{proreg2}.
$\square$
\medskip

Gathering Propositions \ref{proreg1} and \ref{proreg2}, we see that propagation holds for
$L^p$ norms (more precisely, it holds in the spaces $\cup_{q<p, m\in\R}\, L^q_m$) as soon as $p< N/(N-1)$, under the assumptions on $\psi$ in those propositions, and for all dimensions $N$. 
\medskip

We think that propagation in $L^p$ for larger $p$ is
 difficult to treat if one does not want to use the coercivity estimate  
$ \sum_{i,j} (a_{ij} * f) \,\xi_i\,\xi_j \ge C\,(1 +|v|^2)^{\gamma/2}\,|\xi|^2$.
\par
Using this coercivity estimate leads to proofs which are much closer to those of \cite{kcw} than those of Propositions \ref{proreg1} and \ref{proreg2}. We briefly describe how they work: first,
in estimate (\ref{uut}), we use the coercivity estimate and Cauchy-Schwarz inequality in order
to replace
$f^{2k+1}$ by $f^{k+1}$. Then we can take $\beta = 1/(1+k)$, so that a bound in $L^{\infty}(\cup_{\var>0, m\in\R} \, L^{A-\var}_m)$ implies a bound in 
$L^{\infty}(\cup_{\var>0, m\in\R}\, L^{1 + 2\,(A-\var)/N}_m)$, which improves the bound 
of Lemma \ref{ebody2}. In this way, we can get the propagation of 
$L^p$ norms (more precisely, we consider the spaces $\cup_{q<p, m\in\R} L^q_m$) as soon as $p< N/(N-2)$. 
\par
This is sufficient to show that $c * f$ is bounded in $L^{\infty}_{loc}([0,T]\times\R^N)$, so
that using the parabolic form (\ref{parab2}) of the equation, it is possible (for $N=3$) to show that
the bound of $f$ in $L^{\infty}([0,T]; H^1(\R^3)) \cap L^2([0,T]; H^2(\R^3))$
 only depends on the initial norm 
$||f_{in}||_{H^1(\R^3)}$, the initial moments (of any order) and the other parameters
appearing in  Proposition \ref{proreg1}. 
\par
Further smoothness estimates are obtained by differentiating the equation. Note that in all
the steps briefly described above, the dependence with respect to $T$ of the constant is
always polynomial, since no Gronwall-type argument is used.
\medskip

\bigskip

\bigskip

{\bf{Acknowledgment}}: Cl\'ement Mouhot is thanked for fruitful discussions during
the preparation of this work.


\bigskip

\noindent
\begin{thebibliography}{10}

\bibitem{alex_entropy}
R. Alexandre.
\newblock Sur le taux de dissipation d'entropie sans troncature angulaire.
\newblock {\em C.R. Acad. Sci. Paris, S\'erie I}, {\bf{326}}, (1998), 311--315.

\bibitem{alex_new}
R. Alexandre, J. Liao, and C. Lin.
\newblock Some a priori estimates for the homogeneous Landau equation
with soft potentials.
\newblock arXiv:1302.1814.

\bibitem{advw}
R. Alexandre, L. Desvillettes, C. Villani and B. Wennberg.
\newblock Entropy dissipation and long-range interactions.
\newblock {\em Arch. Rat. Mech. Anal.}, {\bf{152}}, (2000), 327-355.
 
 \bibitem{av}
R. Alexandre, and C. Villani.
\newblock On the Landau approximation in plasma physics.
\newblock {\em Annales Inst. Henri Poincar\'e, (C) Analyse non-lin\'eaire},  {\bf{21}} n.1 (2004), 61-95.

\bibitem{arsenev}
A.A. Arsenev and N.V. Peskov.
\newblock On the existence of a generalized solution of Landau's equation.
\newblock {\em Zh. Vychisl. Mat. Mat. Fiz.}, {\bf{17}}, (1977), 1063--1068.

\bibitem{boby}
A. Bobylev, M. Pulvirenti, C. Saffirio.
\newblock From particle systems to the Landau equation: a consistency result.
\newblock {\em Comm. Math. Phys.}, {\bf{319}}, n.3, (2013), 683--702.

\bibitem{bogol}
N.N. Bogolyubov, 
\newblock Problems of a Dynamical Theory in Statistical Physics 
\newblock State Technical Press, 1946 (in Russian). English translation in Studies in Statistical Mechanics I, edited by J. de Boer and G. E. Uhlenbeck, part A, North-Holland, Amsterdam, 1962.

\bibitem{carlencarv:entropy:92}
E.~A.~Carlen and M.~C. Carvalho.
\newblock Strict entropy production bounds and stability of the rate of
  convergence to equilibrium for the {B}oltzmann equation.
\newblock {\em J. Stat. Phys.}, {\bf 67 } n.3-4, (1992), {575--608}.

\bibitem{carlencarv:physic:94}
E.~A. Carlen and M.~C. Carvalho.
\newblock Entropy production estimates for {B}oltzmann equations with
  physically realistic collision kernels.
\newblock {\em J. Stat. Phys.}, {\bf 74} n.3-4, (1994), {743--782}.


\bibitem{chapman}
S.~Chapman and T.G.~Cowling.
\newblock  {\em The mathematical theory of non--uniform gases.}
\newblock Cambridge Univ. Press., London, 1952.

\bibitem{chenlixu1}
H. Chen, W.-X. Li, C.-J. Xu.
\newblock Propagation of Gevrey regularity for solutions of Landau equations,
\newblock {\em Kinetic and Related Models}, {\bf{1}} n.3 (2008), 355--368.

\bibitem{chenlixu2}
H. Chen, W.-X. Li, C.-J. Xu.
\newblock Analytic smoothness effect of solutions for spatially homogeneous Landau equation,
\newblock {\em J. Differ. Equations} {\bf{248}}, (2010), 77--94.

\bibitem{degondlemou:linear:96}
P.~Degond and M.~Lemou.
\newblock Dispersion relations for the linearized {F}okker-{P}lanck equation.
\newblock {\em Arch. Rat. Mech. Anal.}, {\bf{138}}, (1997), 137--167.

\bibitem{desv:entr:89}
L.~Desvillettes.
\newblock Entropy dissipation rate and convergence in kinetic equations.
\newblock {\em Comm. Math. Phys.}, {\bf 123} n.4, (1988), {687--702}.


\bibitem{desv_TTSP} 
L.~Desvillettes.
\newblock On asymptotics of the Boltzmann equation when the collisions become grazing.
\newblock {\em Transport Theory Statist. Phys.}, {\bf{21}}, n.3 (1992), 259--276.


\bibitem{desvvill:I:97}
L.~Desvillettes and C.~Villani.
\newblock On the spatially homogeneous {L}andau equation for hard potentials.
  {P}art {I}. Existence, uniqueness and smoothness.
\newblock {\em Commun. Partial Differential Equations},  {\bf{25}}, n.1-2 (2000), 179-259.


\bibitem{DV2}
L.~Desvillettes and C.~Villani.
\newblock On the spatially homogeneous {L}andau equation for hard potentials.
  {P}art {II}. H-Theorem and applications. 
\newblock {\em Commun. Partial Differential Equations},  {\bf{25}}, n.1-2 (2000), 261-298.


\bibitem{fournier}
N. Fournier.
\newblock Uniqueness of bounded solutions for the homogeneous Landau equation with a Coulomb
potential.
 \newblock {\em Commun. Math. Phys.}, {\bf{299}}, (2010), 765--782.

\bibitem{fournier_guerin}
N. Fournier, H. Gu\'erin.
\newblock Well-posedness of the spatially homogeneous Landau equation for soft potentials.
 \newblock {\em J. Funct. Anal.}, {\bf{25}}, n.8, (2009), 2542--2560.

\bibitem{guerin}
H. Gu\'erin.
\newblock Solving Landau equation for some soft potentials through a probabilistic approach.
\newblock {\em Ann. Appl. Probab.}, {\bf{13}}, n.2, (2003), 515--539.

\bibitem{gs1}
P. T. Gressman and R. M. Strain.
\newblock Global Classical solutions of the Boltzmann equation with Long-Range interactions.
\newblock {\em Proc. Nat. Acad. Sci. U.S.A.}, {\bf{107}}, n.13, (2010), 5744--5749.

\bibitem{gs2}
P. T. Gressman and R. M. Strain.
\newblock Global Classical Solutions of the Boltzmann Equation without Angular Cut-off.
\newblock {\em J. Amer. Math. Soc.}, {\bf{24}}, n.3, (2011), 771--847.

\bibitem{gs3}
P. T. Gressman and R. M. Strain.
\newblock Sharp anisotropic estimates for the Boltzmann collision operator and its entropy production.
\newblock {\em Advances in Math.}, {\bf{227}}, n.6, (2011), 2349--2384.

\bibitem{guo}
Y. Guo.
\newblock The Landau Equation in a Periodic Box.
\newblock {\em  Comm. Math. Phys.}, {\bf{231}}, n.3, (2002), 391-434.


\bibitem{lifschitz}
E.M.~Lifschitz and L.P.~Pitaevskii.
\newblock {\em {Physical kinetics}}.
\newblock Perg. Press., Oxford, 1981.

\bibitem{lions}
P.-L.~Lions.
\newblock Regularity and compactness for Boltzmann collision operators
without angular cut-off. 
\newblock {\em C.R. Acad. Sci. Paris, Serie I}, {\bf{326}}, n.1 (1998), 37--41.


\bibitem{mori}
Y. Morimoto, K. Pravda-Starov and C.-J. Xu.
\newblock A remark on the ultra-analytic smoothing properties of the spatially homogeneous Landau equation.
\newblock {\em Kinet. Relat. Models}, {\bf{6}}, n.4 (2013), 715--727.


\bibitem{toscani_villani}
G.~Toscani and C.~Villani.
\newblock On the trend to equilibrium for some dissipative systems with slowly increasing a priori bounds.
\newblock {\em  J. Statist. Phys.} {\bf{98}}, n.5-6, (2000), 1279--1309.



\bibitem{villani_hb}
C.~Villani.
\newblock A Review of Mathematical Topics in Collisional Kinetic Theory, in Handbook of Mathematical Fluid Dynamics, edited by S. Friedlander and D. Serre, vol. 1, Elsevier, 2002, ISBN 978-0-444-50330-5. doi:10.1016/S1874-5792(02)80004-0.

\bibitem{vill:lm:96}
C.~Villani.
\newblock On the spatially homogeneous {L}andau equation for {M}axwellian
  molecules.
\newblock {\em Math. Meth. Mod. Appl. Sci.} {\bf{8}} n.6, (1998), 957--983.


\bibitem{vill:new:97}
C.~Villani.
\newblock On a new class of weak solutions to the spatially homogeneous
  {B}oltzmann and {L}andau equations.
\newblock {\em Arch. Rat. Mech. Anal.} {\bf{143}} n.3, (1998), 273--307.
 

\bibitem{vill:smooth}
C.~Villani.
\newblock Regularity estimates via the entropy dissipation for the spatially homogeneous Boltzmann equation without cut-off.
\newblock {\em Rev. Matem. Iberoam.} {\bf{15}} n.2, (1999), 335--352. 


\bibitem{kcw}
K.-C. Wu.
\newblock Global in time estimates for the spatially homogeneous Landau equation with soft potentials. 
\newblock {\em J. Funct. Anal.}, {\bf{266}}, (2014), 3134-3155.

\end{thebibliography}
\end{document}